\pdfoutput=1
\documentclass[11pt]{article}
\usepackage[left=1in,right=1in,top=1in,bottom=1in]{geometry}
\usepackage{times}
\usepackage{expl3}
\usepackage{cite}
\usepackage[table]{xcolor}
\usepackage{multirow}
\usepackage{stackengine} 
\usepackage{hhline}
\usepackage{lipsum}
\usepackage{titlesec}
\usepackage{wrapfig}
\usepackage{enumerate}
\usepackage{epsfig}
\usepackage{amsmath}
\usepackage{tabularx}
\usepackage{array}
\usepackage{booktabs}
\usepackage{enumitem}
\usepackage{bbm}
\usepackage{calc}
\usepackage{graphicx}
\usepackage{amsmath}
\usepackage[title]{appendix}
\usepackage{amssymb}
\usepackage{epstopdf}
\usepackage{boldline}
\usepackage{arydshln}
\usepackage{calligra}
\usepackage{bm}
\usepackage{url}
\usepackage{blindtext}
\usepackage{accents}

\newcommand{\define}{\stackrel{\mbox{\tiny def}}{=}}

\newtheorem{definition}{Definition}
\newtheorem{theorem}{Theorem}
\newtheorem{corollary}{Corollary}
\newtheorem{lemma}{Lemma}

\newtheorem{example}{Example}
\newtheorem{remark}{Remark}
\usepackage{mathtools}
\usepackage{epstopdf}
\usepackage{balance}
\usepackage{thmtools}
\usepackage{thm-restate}
\usepackage{hyperref}
\usepackage{cleveref}
\usepackage[mathscr]{euscript}

\usepackage[ruled,vlined]{algorithm2e}
\include{pythonlisting}
\newcommand{\ostar}{\mathbin{\mathpalette\make@circled\star}}

\makeatletter
\newcommand{\removelatexerror}{\let\@latex@error\@gobble}
\makeatother
\makeatletter
\newcommand*{\rom}[1]{\expandafter\@slowromancap\romannumeral #1@}
\makeatother

\ExplSyntaxOn
\newcommand\latinabbrev[1]{
  \peek_meaning:NTF . {
    #1\@}%
  { \peek_catcode:NTF a {
      #1.\@ }%
    {#1.\@}}}
\ExplSyntaxOff

\titleclass{\subsubsubsection}{straight}[\subsubsection]

\begin{document}
\vspace{1cm}
\title{Generalized Converses of Operator Jensen’s Inequalities with Applications to Hypercomplex Function Approximations and Bounds Algebra}\vspace{1.8cm}
\author{Shih-Yu~Chang
\thanks{Shih Yu Chang is with the Department of Applied Data Science,
San Jose State University, San Jose, CA, U. S. A. (e-mail: {\tt
shihyu.chang@sjsu.edu}). 
           }}

\maketitle

\begin{abstract}
Mond and Pecaric proposed a powerful method, namd as MP method, to deal with operator inequalities. However, this method requires a real-valued function to be convex or concave, and the normalized positive linear map between Hilbert spaces. The objective of this study is to extend the MP method by allowing non-convex or non-concave real-valued functions and nonlinear mapping between Hilbert spaces. The Stone–Weierstrass theorem and Kantorovich function are fundamental components employed in generalizing the MP method inequality in this context. Several examples are presented to demonstrate the inequalities obtained from the conventional MP method by requiring convex function with a normalized positive linear map. Various new inequalities regarding hypercomplex functions, i.e., operator-valued functions with operators as arguments, are derived based on the proposed MP method. These inequalities are applied to approximate hypercomplex functions using ratio criteria and difference criteria. Another application of these new inequalities is to establish bounds for hypercomplex functions algebra, i.e., an abelian monoid for the addition or multiplication of hypercomplex functions, and to derive tail bounds for random tensors ensembles addition or multiplication systematically.
\end{abstract}


\section{Introduction}\label{sec: Introduction}

Operator Jensen's inequality extends the concept of conventional Jensen's inequality from real value argument function to operator argument function and is generally used in optimization, especially in convex optimization involving operator variables, e.g., matrices, tensors. The inequality states that for an operator convex function $f$, we have
\begin{eqnarray}\label{eq:operato Jensen}
f(\Phi(\bm{A})) \leq \Phi(f(\bm{A})),
\end{eqnarray}
where $\bm{A}$ is a self-adjoint operator, $\Phi$ is a normalized positive linear map and $\leq$ is by Loewner ordering sense~\cite{furuta2005mond,fujii2012recentMP}. Operator Jensen's inequality is a powerful device in diverse areas, along with machine learning, signal processing, and control principle. In gadget learning, for example, it is used inside the analysis of convex optimization problems related to matrices, consisting of matrix factorization and matrix completion. In signal processing, it finds applications in the design and analysis of convex optimization algorithms for jobs like blind source separation and beamforming. In the field of control science, it is carried out in the evaluation and synthesis of convex control structures regarding matrix variables. Operator Jensen's inequality gives a way to generalize the idea of convex features to operator variables and presents insights into the function convexity related to matrices or tensors, that are frequently encountered in many mathematical and engineering tasks~\cite{chang2024TZF,chang2024generalizedCDJ,chang2023TenLevenberg,chang2023BMSBSVD,chang2023personalized,chang2023TenLeastSquares,
chang2023TenProj,chang2022TLMS,chang2021TenInv,
chang2022TWF,chang2022Kalman,chang2022TCASI,chang2021RLS}.

Mond and Pe{\v{c}}ari{\'c} demonstrated multiple extensions of Kantorovich-type operator inequalities concerning normalized positive linear maps. They highlighted that determining upper bounds for the difference and ratio in Jensen's inequality can be simplified to solving a single-variable maximization or minimization problem by leveraging the concavity of a real-valued function $f$. Building on this approach, they established complementary inequalities to the H{\"o}lder-McCarthy inequality and Kantorovich-type inequalities, provided estimates for the difference and ratio of operator means, and explored various converses of Jensen's inequality applicable to normalized positive linear maps. In cases where $f$ exhibits concavity, they derived the dual problem. This method, known as the Mond-Pe{\v{c}}ari{\'c} (MP) method, has proven to be highly fruitful in the realm of operator inequalities, offering valuable insights and results~\cite{furuta2005mond,fujii2012recentMP}.

The purpose of this work is to extend MP method by allowing $f$ in Eq.~\eqref{eq:operato Jensen} to be non-convex or non-concave and the mapping $\Phi$ in Eq.~\eqref{eq:operato Jensen} to be a nonlinear mapping. The Stone–Weierstrass theorem and Kantorovich function serve as the primary ingredients employed in generalizing MP method inequality in this context. Several examples are shown here to have those inequalities obtained from the conventional MP method by setting $f$ to be a convex function with $\Phi$ as a normalized positive linear map. Various new inequalities about hypercomplex functions, i.e., operator-valued functions with operators as arguments, are derived based on the proposed MP method. These new inequalities are applied to approximate hypercomplex functions by ratio criteria and by difference criteria. The other application of these new inequalities is to establish hypercomplex functions lower/upper bounds algebra, i.e., these bounds form a \emph{abelian monoid} under the addition or the multiplication of hypercomplex functions. Besides, we also can derive random tensors tail bounds for the addition or the multiplication of random tensors ensembles by an uniform way.

The remainder of this paper is organized as follows. In Section~\ref{sec: Generalized Converses of Operator Jensen’s Inequalities with Nonlinear Maps}, fundamental inequalities about hypercomplex functions are established. In Section~\ref{sec: Generalized Converses of Operator Jensen’s Inequalities: Ratio Kind}, generalized converses of operator Jensen’s inequalities for ratio kinds with nonlinear $\Phi$ are established. On the other hand, generalized converses of operator Jensen’s inequalities for different kinds with nonlinear $\Phi$ are established in Section~\ref{sec: Generalized Converses of Operator Jensen’s Inequalities: Difference Kind}. The first application about applying newly derived inequalities to hypercomplex functions approximation is studied in Section~\ref{sec: Hypercomplex Function Approximations}. The second application about applying newly derived inequalities to investigate lower and upper bounds algebra is presented in Section~\ref{sec: Bounds Algebra}.

\textbf{Nomenclature:} 
Inequalities $\geq, >, \leq,$ and $<$, when applied to operators, follow the Loewner ordering. The symbol $\Lambda(\bm{A})$ denotes the spectrum of the operator $\bm{A}$, i.e., the set of eigenvalues of $\bm{A}$. If $\Lambda(\bm{A})$ consists of real numbers, $\min(\Lambda(\bm{A}))$ and $\max(\Lambda(\bm{A}))$ represent the minimum and maximum values within $\Lambda(\bm{A})$, respectively. For given $M>m>0$ and any $r \in \mathbb{R}$ where $r \neq 1$, the Kantorovich function with respect to $m$, $M$, and $r$ is defined as follows:
\begin{eqnarray}\label{eq: Kantorovich function}
\mathscr{K}(m,M,r)&=&\frac{(mM^r - Mm^r)}{(r-1)(M-m)}\left[\frac{(r-1)(M^r-m^r)}{r(mM^r - Mm^r)}\right]^r.
\end{eqnarray}

\section{Fundamental Inequalities for Hypercomplex Functions with Nonlinear Map}\label{sec: Generalized Converses of Operator Jensen’s Inequalities with Nonlinear Maps}

\subsection{Upper Bound and Lower Bound of Function $f$}\label{sec: Upper Bound and Lower Bound of Function f}

Let's start with the Stone-Weierstrass Theorem, which asserts that any continuous real-valued function $f(x)$ defined on the closed interval $[m, M]$, where $m, M \in \mathbb{R}$ and $M > m$, can be uniformly approximated by a polynomial $p(x)$. This approximation ensures that the absolute difference between $f(x)$ and $p(x)$ is less than any given positive value $\epsilon$ across the entire interval $[m, M]$. Mathematically, this difference is bounded by the supremum norm, denoted as $\left\Vert f - p \right\Vert_{\infty}$, which remains less than $\epsilon$~\cite{de1959stone}.

Given a continuous real-valued function $f(x)$ and a positive error bound $\epsilon$, we can employ the Lagrange polynomial interpolation method, grounded in the Stone-Weierstrass Theorem, to ascertain both an upper polynomial $p_{\mathscr{U}}(x) \geq f(x)$ and a lower polynomial $p_{\mathscr{L}}(x) \leq f(x)$ over the interval $[m, M]$. These polynomials are guaranteed to satisfy the following inequalities:
\begin{eqnarray}\label{eq:lower and upper polynomials}
0 &\leq& p_{\mathscr{U}}(x) - f(x)~~\leq~~\epsilon, \nonumber \\
0 &\leq& f(x)-p_{\mathscr{L}}(x)~~\leq~~\epsilon,
\end{eqnarray}
Additionally, in this paper, we assume that $f(\bm{A})$ is a self-adjoint operator if $\bm{A}$ is a self-adjoint operator.

\subsection{Polynomial Map $\Phi$}\label{sec: Polynomial Maps Phi}

Consider two Hilbert spaces $\mathfrak{H}$ and $\mathfrak{K}$. $\mathbb{B}(\mathfrak{H})$ and $\mathbb{B}(\mathfrak{K})$ represent the semi-algebras comprising all bounded linear operators on these respective spaces. In Choi-Davis-Jensen's Inequality, the mapping $\Phi: \mathbb{B}(\mathfrak{H}) \rightarrow \mathbb{B}(\mathfrak{K})$ is defined as a normalized positive linear map. Such a map is precisely defined by  Definition~\ref{def: normalized positive linear map} below~\cite{furuta2005mond,fujii2012recentMP}.
\begin{definition}\label{def: normalized positive linear map}
A map $\Phi: \mathbb{B}(\mathfrak{H}) \rightarrow \mathbb{B}(\mathfrak{K})$ is considered a normalized positive linear map if it satisfies the following conditions:
\begin{itemize}
\item Linearity: $\Phi(a\bm{X}+b\bm{Y})=a\Phi(\bm{X})+b\Phi(\bm{Y})$ for any $a,b \in \mathbb{C}$ and any $\bm{X}, \bm{Y} \in \mathbb{B}(\mathfrak{H})$.
\item Positivity: If $\bm{X}\geq\bm{Y}$, then $\Phi(\bm{X})\geq\Phi(\bm{Y})$.
\item Normalization: $\Phi(\bm{I}_{\mathfrak{H}})=\bm{I}_{\mathfrak{K}}$, where $\bm{I}_{\mathfrak{H}}$ and $\bm{I}_{\mathfrak{K}}$ are the identity operators of the Hilbert spaces $\mathfrak{H}$ and $\mathfrak{K}$, respectively. 
\end{itemize}
\end{definition}

In this study, we will explore a broader class of $\Phi$ by defining $\Phi$ as follows: 
\begin{eqnarray}\label{eq: new phi def}
\Phi(\bm{X}) &=& \bm{V}^{\ast}\left(\sum\limits_{i=0}^{I}a_{i}\bm{X}^{i}\right)\bm{V}\nonumber \\
&=& \bm{V}^{\ast}\left(\sum\limits_{i_{+}\in S_{I_{+}}}a_{i_{+}}\bm{X}^{i_{+}}+\sum\limits_{i_{-}\in S_{I_{-}}}a_{i_{-}}\bm{X}^{i_{-}}\right)\bm{V},
\end{eqnarray}
where, $\bm{V}$ stands as an isometry within $\mathfrak{H}$, adhering to $\bm{V}^{\ast}\bm{V}=\bm{I}_{\mathfrak{H}}$. In $a_{i}$, $a_{i_{+}}$ signifies the nonnegative coefficients, while $a_{i_{-}}$ denotes the negative coefficients. The set of indices corresponding to positive coefficients is denoted as $S_{I_{+}}$, and those for negative coefficients form $S_{I_{-}}$. It's noteworthy that no constraints regarding linearity, positivity, or normalization are imposed on $\Phi$ as outlined in Eq.~\eqref{eq: new phi def}. With this premise, the conventional notion of a normalized positive linear map, delineated in Definition~\ref{def: normalized positive linear map}, emerges as a special case. This is accomplished by configuring the polynomial $\sum\limits_{i=0}^{I}a_{i}\bm{X}^{i}$ as the identity map, wherein all coefficients $a_i$ except $a_1$ are zero.

Contents from Section~\ref{sec: Upper Bound and Lower Bound of Function f} and Section~\ref{sec: Polynomial Maps Phi} are based on~\cite{chang2024generalizedCDJ}, however, we present them here again for self-contained presentation purpose. 

\subsection{Fundamental Inequalities}\label{sec: Fundamental Inequalities}

In this section, general converses of operator Jensen’s inequalities for any polynomial map $\Phi$ will be obtained. Let us recall Lemma 2 from~\cite{chang2024generalizedCDJ}. 
\begin{lemma}\label{lma: phi(f(A)) bounds}
Given a self-adjoint operator $\bm{A}$ with $\Lambda(\bm{A})$, such that 
\begin{eqnarray}\label{eq1: lma: lower and upper bound for f(A)}
0 &\leq& p_{\mathscr{U}}(x) - f(x)~~\leq~~\epsilon, \nonumber \\
0 &\leq& f(x)-p_{\mathscr{L}}(x)~~\leq~~\epsilon,
\end{eqnarray}
for $x \in [\min(\Lambda(\bm{A})), \max(\Lambda(\bm{A}))]$ with polynomials $p_{\mathscr{L}}(x)$ and $p_{\mathscr{U}}(x)$ expressed by 
\begin{eqnarray}\label{eq: lower poly formats}
p_{\mathscr{L}}(x)=\sum\limits_{k=0}^{n_{\mathscr{L}}}\alpha_k x^k, ~~
p_{\mathscr{U}}(x)=\sum\limits_{j=0}^{n_{\mathscr{U}}}\beta_j x^j. 
\end{eqnarray}
Under the definition of $\Phi$ provided by Eq.~\eqref{eq: new phi def}, we have
\begin{eqnarray}\label{eq1:lma: phi(f(A)) bounds}
\Phi(f(\bm{A})) &\leq& \bm{V}^{\ast}\left\{\sum\limits_{i_{+}\in S_{I_{+}}}a_{i_{+}}\mathscr{K}(\min(\Lambda(p_{\mathscr{U}}(\bm{A}))),\max(\Lambda(p_{\mathscr{U}}(\bm{A}))),i_+)p^{i_+}_{\mathscr{U}}(\bm{A})\right. \nonumber \\
&  &\left.+\sum\limits_{i_{-}\in S_{I_{-}}}a_{i_{-}}\mathscr{K}^{-1}(\min(\Lambda(p_{\mathscr{L}}(\bm{A}))),\max(\Lambda(p_{\mathscr{L}}(\bm{A}))),i_-)p^{i_-}_{\mathscr{L}}(\bm{A})\right\}\bm{V}\nonumber \\
&&\define  \bm{V}^{\ast}\mbox{Poly}_{f,\mathscr{U}}(\bm{A})\bm{V},
\end{eqnarray}
where $\Lambda(p_{\mathscr{U}}(\bm{A})$ and $\Lambda(p_{\mathscr{L}}(\bm{A})$ are spectrms of operators  $p_{\mathscr{U}}(\bm{A}$ and $p_{\mathscr{L}}(\bm{A}$, respectively. On the other hand, we also have
\begin{eqnarray}\label{eq2:lma: phi(f(A)) bounds}
\Phi(f(\bm{A})) &\geq& \bm{V}^{\ast}\left\{\sum\limits_{i_{+}\in S_{I_{+}}}a_{i_{+}}\mathscr{K}^{-1}(\min(\Lambda(p_{\mathscr{L}}(\bm{A}))),\max(\Lambda(p_{\mathscr{L}}(\bm{A}))),i_+)p^{i_+}_{\mathscr{L}}(\bm{A})\right. \nonumber \\
&  &\left.+\sum\limits_{i_{-}\in S_{I_{-}}}a_{i_{-}}\mathscr{K}(\min(\Lambda(p_{\mathscr{U}}(\bm{A}))),\max(\Lambda(p_{\mathscr{U}}(\bm{A}))),i_-)p^{i_-}_{\mathscr{U}}(\bm{A})\right\}\bm{V}\nonumber \\
&&\define  \bm{V}^{\ast}\mbox{Poly}_{f,\mathscr{L}}(\bm{A})\bm{V}.
\end{eqnarray}
\end{lemma}
Given $\Lambda(\bm{A}) \in [m,M]$, the value range for the spectrum $\Lambda(\mbox{Poly}_{f,\mathscr{L}}(\bm{A}))$ is represented by $\widetilde{\mbox{Poly}}_{f,\mathscr{L}}(m,M)$. Similarly, the value range for the spectrum $\Lambda(\mbox{Poly}_{f,\mathscr{U}}(\bm{A}))$ is represented by $\widetilde{\mbox{Poly}}_{f,\mathscr{U}}(m,M)$.

The main theorem of this work is presented below. Theorem~\ref{thm: main 2.3} will give operator inequalties, lower and uppber bounds, of functional with respect to $\Phi(f(\bm{A}))$.     
\begin{theorem}\label{thm: main 2.3}
Let $\bm{A}_j$ be self-adjoint operator with $\Lambda(\bm{A}_j) \in [m_j, M_j]$ for real scalars $m_j <  M_j$. The mapping $\Phi: \mathbb{B}(\mathfrak{H}) \rightarrow \mathbb{B}(\mathfrak{K})$ is defined by Eq.~\eqref{eq: new phi def}. The index $j$ is in the range of $1,2,\cdots,k$, and we have a probability vector $\bm{w}=[w_1,w_2,\cdots, w_k]$ with the dimension $k$, i.e., $\sum\limits_{j=1}^{k}w_j = 1$. Let $f$ be any real valued continous functions defined on the range $\bigcup\limits_{j=1}^k [m_j, M_j]$, represented by $f \in \mathcal{C}(\bigcup\limits_{j=1}^k [m_j, M_j])$. Besides, we assume that the function $f$ satisfes thsoe conditions provided by Eq.~\eqref{eq1: lma: lower and upper bound for f(A)} and Eq.~\eqref{eq: lower poly formats}. The function $g$ is also a real valued continous function defined on the range $\left(\bigcup\limits_{j=1}^k w_j\widetilde{\mbox{Poly}}_{f,\mathscr{L}}(m_j,M_j)\right)\bigcup\left(\bigcup\limits_{j=1}^k w_j\widetilde{\mbox{Poly}}_{f,\mathscr{U}}(m_j,M_j)\right)$. We also have a real valued function $F(u,v)$ with operator monotone on the first variable $u$ defined on $U \times V$ such that $f(\bigcup\limits_{j=1}^k [m_j, M_j]) \subset U$, and  $g\left(\left(\bigcup\limits_{j=1}^k w_j\widetilde{\mbox{Poly}}_{f,\mathscr{L}}(m_j,M_j)\right)\bigcup\left(\bigcup\limits_{j=1}^k w_j\widetilde{\mbox{Poly}}_{f,\mathscr{U}}(m_j,M_j)\right)\right) \subset V$. 

Then, we have the following upper bound:
\begin{eqnarray}\label{eq UB: thm:main 2.3}
F\left(\sum\limits_{j=1}^k w_j\Phi(f(\bm{A}_j)),g\left(\sum\limits_{j=1}^k w_j\bm{V}^{\ast}\mbox{Poly}_{f,\mathscr{U}}(\bm{A}_j)\bm{V}\right)\right) \leq
\max\limits_{x\in\bigcup\limits_{j=1}^k w_j\widetilde{\mbox{Poly}}_{f,\mathscr{U}}(m_j,M_j)}F(x,g(x))\bm{1}_{\mathfrak{K}}.
\end{eqnarray}
Similarly, we also have the following lower bound:
\begin{eqnarray}\label{eq LB: thm:main 2.3}
F\left(\sum\limits_{j=1}^k w_j\Phi(f(\bm{A}_j)),g\left(\sum\limits_{j=1}^k w_j\bm{V}^{\ast}\mbox{Poly}_{f,\mathscr{L}}(\bm{A}_j)\bm{V}\right)\right) \geq
\min\limits_{x\in\bigcup\limits_{j=1}^k w_j\widetilde{\mbox{Poly}}_{f,\mathscr{L}}(m_j,M_j)}F(x,g(x))\bm{1}_{\mathfrak{K}}.
\end{eqnarray}
\end{theorem}
\textbf{Proof:}
We begin with the proof for the upper bound provided by Eq.~\eqref{eq UB: thm:main 2.3}. From Lemma~\ref{lma: phi(f(A)) bounds}, we have 
\begin{eqnarray}\label{eq1: thm:main 2.3}
\Phi(f(\bm{A})) \leq \bm{V}^{\ast}\mbox{Poly}_{f,\mathscr{U}}(\bm{A})\bm{V}.
\end{eqnarray}
By replacing $\bm{A}$ with $\bm{A}_j$ in Eq.~\eqref{eq1: thm:main 2.3} and applying $w_j$ with respect to each $\bm{A}_j$, we have 
\begin{eqnarray}\label{eq2: thm:main 2.3}
\sum\limits_{j=1}^k w_j\Phi(f(\bm{A}_j)) \leq \sum\limits_{j=1}^k w_j \bm{V}^{\ast}\mbox{Poly}_{f,\mathscr{U}}(\bm{A}_j)\bm{V}.
\end{eqnarray}
According to the function $F(u,v)$ condition, we have
\begin{eqnarray}\label{eq3: thm:main 2.3}
\lefteqn{F\left(\sum\limits_{j=1}^k w_j\Phi(f(\bm{A}_j)),g\left(\sum\limits_{j=1}^k w_j \bm{V}^{\ast}\mbox{Poly}_{f,\mathscr{U}}(\bm{A}_j)\bm{V}\right)\right)}\nonumber \\
&\leq&F\left(\sum\limits_{j=1}^k w_j \bm{V}^{\ast}\mbox{Poly}_{f,\mathscr{U}}(\bm{A}_j)\bm{V}, g\left(\sum\limits_{j=1}^k w_j \bm{V}^{\ast}\mbox{Poly}_{f,\mathscr{U}}(\bm{A}_j)\bm{V}\right)\right)\nonumber \\
&\leq&\max\limits_{x \in \bigcup\limits_{j=1}^k w_j\widetilde{\mbox{Poly}}_{f,\mathscr{U}}(m_j,M_j)}F(x,g(x))\bm{1}_{\mathfrak{K}},
\end{eqnarray}
where the last inequality comes from that the spectrum $\Lambda\left(\sum\limits_{j=1}^k w_j \bm{V}^{\ast}\mbox{Poly}_{f,\mathscr{U}}(\bm{A}_j)\bm{V}\right)$ is in the range of \\ $\bigcup\limits_{j=1}^k w_j\widetilde{\mbox{Poly}}_{f,\mathscr{U}}(m_j,M_j)$. The desired inequality provided by Eq.~\eqref{eq UB: thm:main 2.3} is established. 

Now, we will prove the lower bound provided by Eq.~\eqref{eq LB: thm:main 2.3}. From Lemma~\ref{lma: phi(f(A)) bounds}, we have 
\begin{eqnarray}\label{eq4: thm:main 2.3}
\Phi(f(\bm{A})) \geq \bm{V}^{\ast}\mbox{Poly}_{f,\mathscr{L}}(\bm{A})\bm{V}.
\end{eqnarray}
By replacing $\bm{A}$ with $\bm{A}_j$ in Eq.~\eqref{eq1: thm:main 2.3} and applying $w_j$ with respect to each $\bm{A}_j$, we have 
\begin{eqnarray}\label{eq5: thm:main 2.3}
\sum\limits_{j=1}^k w_j\Phi(f(\bm{A}_j)) \geq \sum\limits_{j=1}^k w_j \bm{V}^{\ast}\mbox{Poly}_{f,\mathscr{L}}(\bm{A}_j)\bm{V}.
\end{eqnarray}
From the function $F(u,v)$ condition, we also have
\begin{eqnarray}\label{eq6: thm:main 2.3}
\lefteqn{F\left(\sum\limits_{j=1}^k w_j\Phi(f(\bm{A}_j)),g\left(\sum\limits_{j=1}^k w_j \bm{V}^{\ast}\mbox{Poly}_{f,\mathscr{L}}(\bm{A}_j)\bm{V}\right)\right)}\nonumber \\
&\geq&F\left(\sum\limits_{j=1}^k w_j \bm{V}^{\ast}\mbox{Poly}_{f,\mathscr{U}}(\bm{A}_j)\bm{V}, g\left(\sum\limits_{j=1}^k w_j \bm{V}^{\ast}\mbox{Poly}_{f,\mathscr{L}}(\bm{A}_j)\bm{V}\right)\right)\nonumber \\
&\geq&\min\limits_{x\in \bigcup\limits_{j=1}^k w_j\widetilde{\mbox{Poly}}_{f,\mathscr{L}}(m_j,M_j)}F(x,g(x))\bm{1}_{\mathfrak{K}},
\end{eqnarray}
where the last inequality comes from that the spectrum $\Lambda\left(\sum\limits_{j=1}^k w_j \bm{V}^{\ast}\mbox{Poly}_{f,\mathscr{L}}(\bm{A}_j)\bm{V}\right)$ is in the range of \\ $\bigcup\limits_{j=1}^k w_j\widetilde{\mbox{Poly}}_{f,\mathscr{L}}(m_j,M_j)$. The desired inequality provided by Eq.~\eqref{eq LB: thm:main 2.3} is established. 
$\hfill \Box$

By applying the function $F(u,v)$ with the following format:
\begin{eqnarray}\label{eq: F u v }
F(u,v) &=& u - \alpha v,
\end{eqnarray} 
where $\alpha \in \mathbb{R}$, we can have the following Theorem~\ref{thm: 2.4}. Theorem~\ref{thm: 2.4} will provide generalized converses of operator Jensen’s inequalities with ratio and difference kinds in coming sections. 

\begin{theorem}\label{thm: 2.4}
Let $\bm{A}_j$ be self-adjoint operator with $\Lambda(\bm{A}_j) \in [m_j, M_j]$ for real scalars $m_j <  M_j$. The mapping $\Phi: \mathbb{B}(\mathfrak{H}) \rightarrow \mathbb{B}(\mathfrak{K})$ is defined by Eq.~\eqref{eq: new phi def}. The index $j$ is in the range of $1,2,\cdots,k$, and we have a probability vector $\bm{w}=[w_1,w_2,\cdots, w_k]$ with the dimension $k$, i.e., $\sum\limits_{j=1}^{k}w_j = 1$. Let $f$ be any real valued continous functions defined on the range $\bigcup\limits_{j=1}^k [m_j, M_j]$, represented by $f \in \mathcal{C}(\bigcup\limits_{j=1}^k [m_j, M_j])$. Besides, we assume that the function $f$ satisfes thsoe conditions provided by Eq.~\eqref{eq1: lma: lower and upper bound for f(A)} and Eq.~\eqref{eq: lower poly formats}. The function $g$ is also a real valued continous function defined on the range $\left(\bigcup\limits_{j=1}^k w_j\widetilde{\mbox{Poly}}_{f,\mathscr{L}}(m_j,M_j)\right)\bigcup\left(\bigcup\limits_{j=1}^k w_j\widetilde{\mbox{Poly}}_{f,\mathscr{U}}(m_j,M_j)\right)$. We also have a real valued function $F(u,v)$ defined as Eq.~\eqref{eq: F u v } with support domain on $U \times V$ such that $f(\bigcup\limits_{j=1}^k [m_j, M_j]) \subset U$, and  $g\left(\left(\bigcup\limits_{j=1}^k w_j\widetilde{\mbox{Poly}}_{f,\mathscr{L}}(m_j,M_j)\right)\bigcup\left(\bigcup\limits_{j=1}^k w_j\widetilde{\mbox{Poly}}_{f,\mathscr{U}}(m_j,M_j)\right)\right) \subset V$. 

Then, we have the following upper bound:
\begin{eqnarray}\label{eq UB: thm: 2.4}
\sum\limits_{j=1}^k w_j\Phi(f(\bm{A}_j)) \leq \alpha g\left(\sum\limits_{j=1}^k w_j\bm{V}^{\ast}\mbox{Poly}_{f,\mathscr{U}}(\bm{A}_j)\bm{V}\right)+
\max\limits_{x\in\bigcup\limits_{j=1}^k w_j\widetilde{\mbox{Poly}}_{f,\mathscr{U}}(m_j,M_j)}(x-\alpha g(x))\bm{1}_{\mathfrak{K}}.
\end{eqnarray}
Similarly, we also have the following lower bound:
\begin{eqnarray}\label{eq LB: thm: 2.4}
\sum\limits_{j=1}^k w_j\Phi(f(\bm{A}_j)) \geq \alpha g\left(\sum\limits_{j=1}^k w_j\bm{V}^{\ast}\mbox{Poly}_{f,\mathscr{L}}(\bm{A}_j)\bm{V}\right) +
\min\limits_{x\in\bigcup\limits_{j=1}^k w_j\widetilde{\mbox{Poly}}_{f,\mathscr{L}}(m_j,M_j)}(x-\alpha g(x))\bm{1}_{\mathfrak{K}}.
\end{eqnarray}
\end{theorem} 
\textbf{Proof:}
By setting $F(u,v)=u - \alpha v$ in Eq.~\eqref{eq UB: thm:main 2.3} in Theorem~\ref{thm: main 2.3}, we have the desired inequality provided by Eq.~\eqref{eq UB: thm: 2.4}. Similarly, By setting $F(u,v)=u - \alpha v$ in Eq.~\eqref{eq LB: thm:main 2.3} in Theorem~\ref{thm: main 2.3}, we have the desired inequality provided by Eq.~\eqref{eq LB: thm: 2.4}.
$\hfill \Box$

Corollary~\ref{cor: special cases of g func} below is provided to give upper and lower bounds for special types of the function $g$ by applying Theorem~\ref{thm: 2.4}.

\begin{corollary}\label{cor: special cases of g func}
Let $\bm{A}_j$ be self-adjoint operator with $\Lambda(\bm{A}_j) \in [m_j, M_j]$ for real scalars $m_j <  M_j$. The mapping $\Phi: \mathbb{B}(\mathfrak{H}) \rightarrow \mathbb{B}(\mathfrak{K})$ is defined by Eq.~\eqref{eq: new phi def}. The index $j$ is in the range of $1,2,\cdots,k$, and we have a probability vector $\bm{w}=[w_1,w_2,\cdots, w_k]$ with the dimension $k$, i.e., $\sum\limits_{j=1}^{k}w_j = 1$. Let $f$ be any real valued continous functions defined on the range $\bigcup\limits_{j=1}^k [m_j, M_j]$, represented by $f \in \mathcal{C}(\bigcup\limits_{j=1}^k [m_j, M_j])$. Besides, we assume that the function $f$ satisfes thsoe conditions provided by Eq.~\eqref{eq1: lma: lower and upper bound for f(A)} and Eq.~\eqref{eq: lower poly formats}. 

(I) If $g(x) = x^q$, where $q \in \mathbb{R}$ and $\left(\left(\bigcup\limits_{j=1}^k w_j\widetilde{\mbox{Poly}}_{f,\mathscr{L}}(m_j,M_j)\right)\bigcup\left(\bigcup\limits_{j=1}^k w_j\widetilde{\mbox{Poly}}_{f,\mathscr{U}}(m_j,M_j)\right)\right) \geq 0$, we have the upper bound for $\sum\limits_{j=1}^k w_j\Phi(f(\bm{A}_j))$:
\begin{eqnarray}\label{eq power U: cor: special cases of g func}
\sum\limits_{j=1}^k w_j\Phi(f(\bm{A}_j)) \leq \alpha \left(\sum\limits_{j=1}^k w_j\bm{V}^{\ast}\mbox{Poly}_{f,\mathscr{L}}(\bm{A}_j)\bm{V}\right)^q +
\max\limits_{x\in\bigcup\limits_{j=1}^k w_j\widetilde{\mbox{Poly}}_{f,\mathscr{U}}(m_j,M_j)}(x-\alpha x^q)\bm{1}_{\mathfrak{K}},
\end{eqnarray}
and we have the lower bound for $\sum\limits_{j=1}^k w_j\Phi(f(\bm{A}_j))$:
\begin{eqnarray}\label{eq power L: cor: special cases of g func}
\sum\limits_{j=1}^k w_j\Phi(f(\bm{A}_j)) \geq \alpha \left(\sum\limits_{j=1}^k w_j\bm{V}^{\ast}\mbox{Poly}_{f,\mathscr{L}}(\bm{A}_j)\bm{V}\right)^q +
\min\limits_{x\in\bigcup\limits_{j=1}^k w_j\widetilde{\mbox{Poly}}_{f,\mathscr{L}}(m_j,M_j)}(x-\alpha x^q)\bm{1}_{\mathfrak{K}}.
\end{eqnarray}

(II) If $g(x) = \log(x)$ and $\left(\left(\bigcup\limits_{j=1}^k w_j\widetilde{\mbox{Poly}}_{f,\mathscr{L}}(m_j,M_j)\right)\bigcup\left(\bigcup\limits_{j=1}^k w_j\widetilde{\mbox{Poly}}_{f,\mathscr{U}}(m_j,M_j)\right)\right) > 0$ , we have the upper bound for $\sum\limits_{j=1}^k w_j\Phi(f(\bm{A}_j))$:
\begin{eqnarray}\label{eq log U: cor: special cases of g func}
\sum\limits_{j=1}^k w_j\Phi(f(\bm{A}_j)) \leq \alpha \log\left(\sum\limits_{j=1}^k w_j\bm{V}^{\ast}\mbox{Poly}_{f,\mathscr{L}}(\bm{A}_j)\bm{V}\right) +
\max\limits_{x\in\bigcup\limits_{j=1}^k w_j\widetilde{\mbox{Poly}}_{f,\mathscr{U}}(m_j,M_j)}(x-\alpha \log(x)\bm{1}_{\mathfrak{K}},
\end{eqnarray}
and we have the lower bound for $\sum\limits_{j=1}^k w_j\Phi(f(\bm{A}_j))$:
\begin{eqnarray}\label{eq log L: cor: special cases of g func}
\sum\limits_{j=1}^k w_j\Phi(f(\bm{A}_j)) \geq \alpha \log\left(\sum\limits_{j=1}^k w_j\bm{V}^{\ast}\mbox{Poly}_{f,\mathscr{L}}(\bm{A}_j)\bm{V}\right) +
\min\limits_{x\in\bigcup\limits_{j=1}^k w_j\widetilde{\mbox{Poly}}_{f,\mathscr{L}}(m_j,M_j)}(x-\alpha \log(x))\bm{1}_{\mathfrak{K}}.
\end{eqnarray}

(III) If $g(x) = \exp(x)$ , we have the upper bound for $\sum\limits_{j=1}^k w_j\Phi(f(\bm{A}_j))$:
\begin{eqnarray}\label{eq exp U: cor: special cases of g func}
\sum\limits_{j=1}^k w_j\Phi(f(\bm{A}_j)) \leq \alpha \exp\left(\sum\limits_{j=1 }^k w_j\bm{V}^{\ast}\mbox{Poly}_{f,\mathscr{L}}(\bm{A}_j)\bm{V}\right) +
\max\limits_{x\in\bigcup\limits_{j=1}^k w_j\widetilde{\mbox{Poly}}_{f,\mathscr{U}}(m_j,M_j)}(x-\alpha \exp(x))\bm{1}_{\mathfrak{K}},
\end{eqnarray}
and we have the lower bound for $\sum\limits_{j=1}^k w_j\Phi(f(\bm{A}_j))$:
\begin{eqnarray}\label{eq exp L: cor: special cases of g func}
\sum\limits_{j=1}^k w_j\Phi(f(\bm{A}_j)) \geq \alpha \exp\left(\sum\limits_{j=1}^k w_j\bm{V}^{\ast}\mbox{Poly}_{f,\mathscr{L}}(\bm{A}_j)\bm{V}\right) +
\min\limits_{x\in\bigcup\limits_{j=1}^k w_j\widetilde{\mbox{Poly}}_{f,\mathscr{L}}(m_j,M_j)}(x-\alpha \exp(x))\bm{1}_{\mathfrak{K}}.
\end{eqnarray}
\end{corollary}

In the following Example~\ref{exp: bound by linear with positive a_i}, we will assume that the function $f$
is bounded by linear functions and derive related inequalties given by Theorem~\ref{thm: main 2.3}, Theorem~\ref{thm: 2.4}, and Corollary~\ref{cor: special cases of g func}.

\begin{example}\label{exp: bound by linear with positive a_i}
In this example, we assume that the function $f$ satsifies the following:
\begin{eqnarray}\label{eq1: exp: bound by linear with positive a_i}
0 &\leq&\overbrace{(\beta_0 + \beta_1 x)}^{p_{\mathscr{U}}(x)} - f(x)~~\leq~~\epsilon, \nonumber \\
0 &\leq& f(x)-\overbrace{(\alpha_0 + \alpha_1 x)}^{p_{\mathscr{L}}(x)}~~\leq~~\epsilon,
\end{eqnarray}
where $x \in [m,M]$ with positive $\alpha_0, \alpha_1, \beta_0$ and $\beta_1$. We also assume that $\Lambda(\bm{A}_j) \in [m,M]$ and $f(\bm{A}_j) \geq \bm{0}$ for all $j \in 1,2,\cdots,k$. Moreover, the mapping $\Phi$ defined by Eq.~\eqref{eq: new phi def} has all coefficients $a_i \geq 0$. Then, from Lemma~\ref{lma: phi(f(A)) bounds}, we have 
\begin{eqnarray}\label{eq2 UB: exp: bound by linear with positive a_i}
\Phi(f(\bm{A})) &\leq& \bm{V}^{\ast}\left\{\sum\limits_{i=0}^{I}a_{i}\mathscr{K}(\beta_0 + \beta_1 m,\beta_0 + \beta_1 M, i)(\beta_0 + \beta_1 \bm{A})^i \right\}\bm{V}.
\end{eqnarray}
Besides, we also have
\begin{eqnarray}\label{eq2 LB: exp: bound by linear with positive a_i}
\Phi(f(\bm{A})) &\geq& \bm{V}^{\ast}\left\{\sum\limits_{i=0}^{I}a_{i}\mathscr{K}^{-1}(\alpha_0 + \alpha_1 m,\alpha_0 + \alpha_1 M, i)(\alpha_0 + \alpha_1 \bm{A})^i \right\}\bm{V}.
\end{eqnarray}
From Eq.~\eqref{eq2 UB: exp: bound by linear with positive a_i}, we can express $\mbox{Poly}_{f,\mathscr{U}}(\bm{A})$ as:
\begin{eqnarray}\label{eq3 UB: exp: bound by linear with positive a_i}
\mbox{Poly}_{f,\mathscr{U}}(\bm{A})&=& \sum\limits_{i=0}^I c_i \bm{A}^i,
\end{eqnarray}
where coeffcients  $c_i$ are given by 
\begin{eqnarray}\label{eq3-1 UB: exp: bound by linear with positive a_i}
c_i&=&\sum\limits_{j=0}^{I-i}{i+j \choose i}a_{i+j} \mathscr{K}(\beta_0 + \beta_1 m,\beta_0 + \beta_1 M, i+j) \beta^i_1 \beta^j_0.
\end{eqnarray} 
Similarly, from Eq.~\eqref{eq2 LB: exp: bound by linear with positive a_i}, we can express $\mbox{Poly}_{f,\mathscr{L}}(\bm{A})$ as:
\begin{eqnarray}\label{eq3 LB: exp: bound by linear with positive a_i}
\mbox{Poly}_{f,\mathscr{L}}(\bm{A})&=& \sum\limits_{i=0}^I d_i \bm{A}^i,
\end{eqnarray}
where coeffcients  $c_i$ are given by 
\begin{eqnarray}\label{eq3-1 LB: exp: bound by linear with positive a_i}
d_i&=&\sum\limits_{j=0}^{I-i}{i+j \choose i}a_{i+j} \mathscr{K}^{-1}(\alpha_0 + \alpha_1 m,\alpha_0 + \alpha_1 M, i+j) \alpha^i_1 \alpha^j_0.
\end{eqnarray} 

Suppose conditions provided by Theorem~\ref{thm: main 2.3} are valid, from polynomials given by Eq.~\eqref{eq3 UB: exp: bound by linear with positive a_i} and Eq.~\eqref{eq3 LB: exp: bound by linear with positive a_i}, we have 
\begin{eqnarray}\label{eq UB: thm:main 2.3 exp1}
F\left(\sum\limits_{j=1}^k w_j\Phi(f(\bm{A}_j)),g\left(\sum\limits_{j=1}^k w_j\bm{V}^{\ast}\mbox{Poly}_{f,\mathscr{U}}(\bm{A}_j)\bm{V}\right)\right) \leq
\max\limits_{x\in\widetilde{\mbox{Poly}}_{f,\mathscr{U}}(m,M)}F(x,g(x))\bm{1}_{\mathfrak{K}}.
\end{eqnarray}
Similarly, we also have the following lower bound:
\begin{eqnarray}\label{eq LB: thm:main 2.3 exp1}
F\left(\sum\limits_{j=1}^k w_j\Phi(f(\bm{A}_j)),g\left(\sum\limits_{j=1}^k w_j\bm{V}^{\ast}\mbox{Poly}_{f,\mathscr{L}}(\bm{A}_j)\bm{V}\right)\right) \geq
\min\limits_{x\in\widetilde{\mbox{Poly}}_{f,\mathscr{L}}(m,M)}F(x,g(x))\bm{1}_{\mathfrak{K}}.
\end{eqnarray}
If $F(u,v) = u - \alpha v$, we can apply Theorem~\ref{thm: 2.4} to the function $f$ and $\Phi$ provided in this example to obtain the following:
\begin{eqnarray}\label{eq UB: thm: 2.4 exp1}
\sum\limits_{j=1}^k w_j\Phi(f(\bm{A}_j)) \leq \alpha g\left(\sum\limits_{j=1}^k w_j\bm{V}^{\ast}\mbox{Poly}_{f,\mathscr{U}}(\bm{A}_j)\bm{V}\right)+
\max\limits_{x\in\widetilde{\mbox{Poly}}_{f,\mathscr{U}}(m,M)}(x-\alpha g(x))\bm{1}_{\mathfrak{K}}.
\end{eqnarray}
Similarly, we also have the following lower bound:
\begin{eqnarray}\label{eq LB: thm: 2.4 exp1}
\sum\limits_{j=1}^k w_j\Phi(f(\bm{A}_j)) \geq \alpha g\left(\sum\limits_{j=1}^k w_j\bm{V}^{\ast}\mbox{Poly}_{f,\mathscr{L}}(\bm{A}_j)\bm{V}\right) +
\min\limits_{x\in\widetilde{\mbox{Poly}}_{f,\mathscr{L}}(m,M)}(x-\alpha g(x))\bm{1}_{\mathfrak{K}}.
\end{eqnarray}
Finally, the application of Corollary~\ref{cor: special cases of g func} to the function $f$ and $\Phi$ provided in this example will get:

(I) If $g(x) = x^q$, where $q \in \mathbb{R}$ and $\left(\widetilde{\mbox{Poly}}_{f,\mathscr{L}}(m,M)\bigcup\widetilde{\mbox{Poly}}_{f,\mathscr{U}}(m,M)\right) \geq 0$, we have the upper bound for $\sum\limits_{j=1}^k w_j\Phi(f(\bm{A}_j))$:
\begin{eqnarray}\label{eq power U: cor: special cases of g func exp1}
\sum\limits_{j=1}^k w_j\Phi(f(\bm{A}_j)) \leq \alpha \left(\sum\limits_{j=1}^k w_j\bm{V}^{\ast}\mbox{Poly}_{f,\mathscr{U}}(\bm{A}_j)\bm{V}\right)^q +
\max\limits_{x\in\widetilde{\mbox{Poly}}_{f,\mathscr{U}}(m,M)}(x-\alpha x^q)\bm{1}_{\mathfrak{K}},
\end{eqnarray}
and we have the lower bound for $\sum\limits_{j=1}^k w_j\Phi(f(\bm{A}_j))$:
\begin{eqnarray}\label{eq power L: cor: special cases of g func exp1}
\sum\limits_{j=1}^k w_j\Phi(f(\bm{A}_j)) \geq \alpha \left(\sum\limits_{j=1}^k w_j\bm{V}^{\ast}\mbox{Poly}_{f,\mathscr{L}}(\bm{A}_j)\bm{V}\right)^q +
\min\limits_{x\in\widetilde{\mbox{Poly}}_{f,\mathscr{L}}(m,M)}(x-\alpha x^q)\bm{1}_{\mathfrak{K}}.
\end{eqnarray}

(II) If $g(x) = \log(x)$ and $\left(\widetilde{\mbox{Poly}}_{f,\mathscr{L}}(m,M)\bigcup\widetilde{\mbox{Poly}}_{f,\mathscr{U}}(m,M)\right) > 0$ , we have the upper bound for $\sum\limits_{j=1}^k w_j\Phi(f(\bm{A}_j))$:
\begin{eqnarray}\label{eq log U: cor: special cases of g func exp1}
\sum\limits_{j=1}^k w_j\Phi(f(\bm{A}_j)) \leq \alpha \log\left(\sum\limits_{j=1}^k w_j\bm{V}^{\ast}\mbox{Poly}_{f,\mathscr{U}}(\bm{A}_j)\bm{V}\right) +
\max\limits_{x\in\widetilde{\mbox{Poly}}_{f,\mathscr{U}}(m,M)}(x-\alpha \log(x)\bm{1}_{\mathfrak{K}},
\end{eqnarray}
and we have the lower bound for $\sum\limits_{j=1}^k w_j\Phi(f(\bm{A}_j))$:
\begin{eqnarray}\label{eq log L: cor: special cases of g func exp1}
\sum\limits_{j=1}^k w_j\Phi(f(\bm{A}_j)) \geq \alpha \log\left(\sum\limits_{j=1}^k w_j\bm{V}^{\ast}\mbox{Poly}_{f,\mathscr{L}}(\bm{A}_j)\bm{V}\right) +
\min\limits_{x\in\widetilde{\mbox{Poly}}_{f,\mathscr{L}}(m,M)}(x-\alpha \log(x))\bm{1}_{\mathfrak{K}}.
\end{eqnarray}

(III) If $g(x) = \exp(x)$ , we have the upper bound for $\sum\limits_{j=1}^k w_j\Phi(f(\bm{A}_j))$:
\begin{eqnarray}\label{eq exp U: cor: special cases of g func exp1}
\sum\limits_{j=1}^k w_j\Phi(f(\bm{A}_j)) \leq \alpha \exp\left(\sum\limits_{j=1 }^k w_j\bm{V}^{\ast}\mbox{Poly}_{f,\mathscr{U}}(\bm{A}_j)\bm{V}\right) +
\max\limits_{x\in\widetilde{\mbox{Poly}}_{f,\mathscr{U}}(m,M)}(x-\alpha \exp(x))\bm{1}_{\mathfrak{K}},
\end{eqnarray}
and we have the lower bound for $\sum\limits_{j=1}^k w_j\Phi(f(\bm{A}_j))$:
\begin{eqnarray}\label{eq exp L: cor: special cases of g func exp1}
\sum\limits_{j=1}^k w_j\Phi(f(\bm{A}_j)) \geq \alpha \exp\left(\sum\limits_{j=1}^k w_j\bm{V}^{\ast}\mbox{Poly}_{f,\mathscr{L}}(\bm{A}_j)\bm{V}\right) +
\min\limits_{x\in\widetilde{\mbox{Poly}}_{f,\mathscr{L}}(m,M)}(x-\alpha \exp(x))\bm{1}_{\mathfrak{K}}.
\end{eqnarray}

\end{example}

\section{Generalized Converses of Operator Jensen’s Inequalities: Ratio Kind}\label{sec: Generalized Converses of Operator Jensen’s Inequalities: Ratio Kind}

In this section, we will derive the lower and upper bounds for $\sum\limits_{j=1}^k w_j\Phi(f(\bm{A}_j))$ in terms of ratio criteria related to the function $g$.
\begin{theorem}\label{thm: 2.9}
Let $\bm{A}_j$ be self-adjoint operator with $\Lambda(\bm{A}_j) \in [m_j, M_j]$ for real scalars $m_j <  M_j$. The mapping $\Phi: \mathbb{B}(\mathfrak{H}) \rightarrow \mathbb{B}(\mathfrak{K})$ is defined by Eq.~\eqref{eq: new phi def}. The index $j$ is in the range of $1,2,\cdots,k$, and we have a probability vector $\bm{w}=[w_1,w_2,\cdots, w_k]$ with the dimension $k$, i.e., $\sum\limits_{j=1}^{k}w_j = 1$. Let $f$ be any real valued continous functions defined on the range $\bigcup\limits_{j=1}^k [m_j, M_j]$, represented by $f \in \mathcal{C}(\bigcup\limits_{j=1}^k [m_j, M_j])$. Besides, we assume that the function $f$ satisfes thsoe conditions provided by Eq.~\eqref{eq1: lma: lower and upper bound for f(A)} and Eq.~\eqref{eq: lower poly formats}. The function $g$ is also a real valued continous function defined on the range $\left(\bigcup\limits_{j=1}^k w_j\widetilde{\mbox{Poly}}_{f,\mathscr{L}}(m_j,M_j)\right)\bigcup\left(\bigcup\limits_{j=1}^k w_j\widetilde{\mbox{Poly}}_{f,\mathscr{U}}(m_j,M_j)\right)$ and $g(x)\neq 0$ for $x \in \left(\bigcup\limits_{j=1}^k w_j\widetilde{\mbox{Poly}}_{f,\mathscr{L}}(m_j,M_j)\right)\bigcup\left(\bigcup\limits_{j=1}^k w_j\widetilde{\mbox{Poly}}_{f,\mathscr{U}}(m_j,M_j)\right)$. 

(I) If we also assume that $g\left(\sum\limits_{j=1}^k w_j\bm{V}^{\ast}\mbox{Poly}_{f,\mathscr{U}}(\bm{A}_j)\bm{V}\right) > \bm{O}$ and $g\left(\sum\limits_{j=1}^k w_j\bm{V}^{\ast}\mbox{Poly}_{f,\mathscr{L}}(\bm{A}_j)\bm{V}\right) > \bm{O}$, then, we have the following upper bound for $\sum\limits_{j=1}^k w_j\Phi(f(\bm{A}_j))$:
\begin{eqnarray}\label{eq1 UB: thm: 2.9 pos g}
\sum\limits_{j=1}^k w_j\Phi(f(\bm{A}_j))&\leq&\left[\max\limits_{x\in\bigcup\limits_{j=1}^k w_j\widetilde{\mbox{Poly}}_{f,\mathscr{U}}(m_j,M_j)}xg^{-1}(x)\right]g\left(\sum\limits_{j=1}^k w_j\bm{V}^{\ast}\mbox{Poly}_{f,\mathscr{U}}(\bm{A}_j)\bm{V}\right);
\end{eqnarray}
and, the following lower bound for $\sum\limits_{j=1}^k w_j\Phi(f(\bm{A}_j))$:
\begin{eqnarray}\label{eq1 LB: thm: 2.9 pos g}
\sum\limits_{j=1}^k w_j\Phi(f(\bm{A}_j))&\geq&\left[\min\limits_{x\in\bigcup\limits_{j=1}^k w_j\widetilde{\mbox{Poly}}_{f,\mathscr{L}}(m_j,M_j)}xg^{-1}(x)\right]g\left(\sum\limits_{j=1}^k w_j\bm{V}^{\ast}\mbox{Poly}_{f,\mathscr{L}}(\bm{A}_j)\bm{V}\right).
\end{eqnarray}

(II) If we also assume that $g\left(\sum\limits_{j=1}^k w_j\bm{V}^{\ast}\mbox{Poly}_{f,\mathscr{U}}(\bm{A}_j)\bm{V}\right) < \bm{O}$ and $g\left(\sum\limits_{j=1}^k w_j\bm{V}^{\ast}\mbox{Poly}_{f,\mathscr{L}}(\bm{A}_j)\bm{V}\right) < \bm{O}$, then, we have the following upper bound for $\sum\limits_{j=1}^k w_j\Phi(f(\bm{A}_j))$:
\begin{eqnarray}\label{eq1 UB: thm: 2.9 neg g}
\sum\limits_{j=1}^k w_j\Phi(f(\bm{A}_j))&\leq&\left[\min\limits_{x\in\bigcup\limits_{j=1}^k w_j\widetilde{\mbox{Poly}}_{f,\mathscr{L}}(m_j,M_j)}xg^{-1}(x)\right]g\left(\sum\limits_{j=1}^k w_j\bm{V}^{\ast}\mbox{Poly}_{f,\mathscr{L}}(\bm{A}_j)\bm{V}\right);
\end{eqnarray}
and, the following lower bound for $\sum\limits_{j=1}^k w_j\Phi(f(\bm{A}_j))$:
\begin{eqnarray}\label{eq1 LB: thm: 2.9 neg g}
\sum\limits_{j=1}^k w_j\Phi(f(\bm{A}_j))&\geq&\left[\max\limits_{x\in\bigcup\limits_{j=1}^k w_j\widetilde{\mbox{Poly}}_{f,\mathscr{U}}(m_j,M_j)}xg^{-1}(x)\right]g\left(\sum\limits_{j=1}^k w_j\bm{V}^{\ast}\mbox{Poly}_{f,\mathscr{U}}(\bm{A}_j)\bm{V}\right).
\end{eqnarray}
\end{theorem}
\textbf{Proof:}
For part (I), we will apply $F(u,v)$ as 
\begin{eqnarray}\label{eq: F u v 2}
F(u,v)&=&v^{-1/2} u v^{-1/2},
\end{eqnarray}
to Eq.~\eqref{eq UB: thm:main 2.3} in Theorem~\ref{thm: main 2.3}, then, we will obtain 
\begin{eqnarray}\label{eq2 UB: thm: 2.9 pos g}
\left(g\left(\sum\limits_{j=1}^k w_j\bm{V}^{\ast}\mbox{Poly}_{f,\mathscr{U}}(\bm{A}_j)\bm{V}\right)\right)^{-1/2}
\left(\sum\limits_{j=1}^k w_j\Phi(f(\bm{A}_j))\right)
\left(g\left(\sum\limits_{j=1}^k w_j\bm{V}^{\ast}\mbox{Poly}_{f,\mathscr{U}}(\bm{A}_j)\bm{V}\right)\right)^{-1/2}
\nonumber \\
\leq\left[\max\limits_{x\in\bigcup\limits_{j=1}^k w_j\widetilde{\mbox{Poly}}_{f,\mathscr{U}}(m_j,M_j)}xg^{-1}(x)\right]\bm{1}_{\mathfrak{K}}.
\end{eqnarray}
By multiplying $\left(g\left(\sum\limits_{j=1}^k w_j\bm{V}^{\ast}\mbox{Poly}_{f,\mathscr{U}}(\bm{A}_j)\bm{V}\right)\right)^{1/2}$ at both sides of Eq.~\eqref{eq2 UB: thm: 2.9 pos g}, we obtain the desired inequality provided by Eq.~\eqref{eq1 UB: thm: 2.9 pos g}. By applying $F(u,v)$ with Eq.~\eqref{eq: F u v 2} again to Eq.~\eqref{eq LB: thm:main 2.3} in Theorem~\ref{thm: main 2.3}, then, we will obtain 
\begin{eqnarray}\label{eq2 LB: thm: 2.9 pos g}
\left(g\left(\sum\limits_{j=1}^k w_j\bm{V}^{\ast}\mbox{Poly}_{f,\mathscr{L}}(\bm{A}_j)\bm{V}\right)\right)^{-1/2}
\left(\sum\limits_{j=1}^k w_j\Phi(f(\bm{A}_j))\right)
\left(g\left(\sum\limits_{j=1}^k w_j\bm{V}^{\ast}\mbox{Poly}_{f,\mathscr{L}}(\bm{A}_j)\bm{V}\right)\right)^{-1/2}
\nonumber \\
\geq\left[\min\limits_{x\in\bigcup\limits_{j=1}^k w_j\widetilde{\mbox{Poly}}_{f,\mathscr{L}}(m_j,M_j)}xg^{-1}(x)\right]\bm{1}_{\mathfrak{K}}.
\end{eqnarray}
By multiplying $\left(g\left(\sum\limits_{j=1}^k w_j\bm{V}^{\ast}\mbox{Poly}_{f,\mathscr{L}}(\bm{A}_j)\bm{V}\right)\right)^{-1/2}$ at both sides of Eq.~\eqref{eq2 LB: thm: 2.9 pos g}, we obtain the desired inequality provided by Eq.~\eqref{eq1 LB: thm: 2.9 pos g}. 

The proof of Part (II) is immediate obtained by setting $g(x)$ as $-g(x)$ in Eq.~\eqref{eq2 UB: thm: 2.9 pos g} and Eq.~\eqref{eq2 LB: thm: 2.9 pos g}. 
$\hfill \Box$

Next Corollary~\ref{cor: 2.11-14} is obtained by applying Theorem~\ref{thm: 2.9} to special types of function $g$. 
\begin{corollary}\label{cor: 2.11-14}
Let $\bm{A}_j$ be self-adjoint operator with $\Lambda(\bm{A}_j) \in [m_j, M_j]$ for real scalars $m_j <  M_j$. The mapping $\Phi: \mathbb{B}(\mathfrak{H}) \rightarrow \mathbb{B}(\mathfrak{K})$ is defined by Eq.~\eqref{eq: new phi def}. The index $j$ is in the range of $1,2,\cdots,k$, and we have a probability vector $\bm{w}=[w_1,w_2,\cdots, w_k]$ with the dimension $k$, i.e., $\sum\limits_{j=1}^{k}w_j = 1$. Let $f$ be any real valued continous functions defined on the range $\bigcup\limits_{j=1}^k [m_j, M_j]$, represented by $f \in \mathcal{C}(\bigcup\limits_{j=1}^k [m_j, M_j])$. Besides, we assume that the function $f$ satisfes thsoe conditions provided by Eq.~\eqref{eq1: lma: lower and upper bound for f(A)} and Eq.~\eqref{eq: lower poly formats}.

(I) If we also assume that $\left(\sum\limits_{j=1}^k w_j\bm{V}^{\ast}\mbox{Poly}_{f,\mathscr{U}}(\bm{A}_j)\bm{V}\right)^q > \bm{O}$ and $\left(\sum\limits_{j=1}^k w_j\bm{V}^{\ast}\mbox{Poly}_{f,\mathscr{L}}(\bm{A}_j)\bm{V}\right)^q > \bm{O}$ for $q \in \mathbb{R}$, then, we have the following upper bound for $\sum\limits_{j=1}^k w_j\Phi(f(\bm{A}_j))$:
\begin{eqnarray}\label{eq x q UB: cor: 2.11-14}
\sum\limits_{j=1}^k w_j\Phi(f(\bm{A}_j))&\leq&\left[\max\limits_{x\in\bigcup\limits_{j=1}^k w_j\widetilde{\mbox{Poly}}_{f,\mathscr{U}}(m_j,M_j)}x^{1-q}\right]\left(\sum\limits_{j=1}^k w_j\bm{V}^{\ast}\mbox{Poly}_{f,\mathscr{U}}(\bm{A}_j)\bm{V}\right)^q;
\end{eqnarray}
and, the following lower bound for $\sum\limits_{j=1}^k w_j\Phi(f(\bm{A}_j))$:
\begin{eqnarray}\label{eq x q LB: cor: 2.11-14}
\sum\limits_{j=1}^k w_j\Phi(f(\bm{A}_j))&\geq&\left[\min\limits_{x\in\bigcup\limits_{j=1}^k w_j\widetilde{\mbox{Poly}}_{f,\mathscr{L}}(m_j,M_j)}x^{1-q}\right]\left(\sum\limits_{j=1}^k w_j\bm{V}^{\ast}\mbox{Poly}_{f,\mathscr{L}}(\bm{A}_j)\bm{V}\right)^q.
\end{eqnarray}

(II) If we also assume that $\log\left(\sum\limits_{j=1}^k w_j\bm{V}^{\ast}\mbox{Poly}_{f,\mathscr{U}}(\bm{A}_j)\bm{V}\right) > \bm{O}$ and $\log\left(\sum\limits_{j=1}^k w_j\bm{V}^{\ast}\mbox{Poly}_{f,\mathscr{L}}(\bm{A}_j)\bm{V}\right) > \bm{O}$ for $q \in \mathbb{R}$, then, we have the following upper bound for $\sum\limits_{j=1}^k w_j\Phi(f(\bm{A}_j))$:
\begin{eqnarray}\label{eq log x l0 UB: cor: 2.11-14}
\sum\limits_{j=1}^k w_j\Phi(f(\bm{A}_j))&\leq&\left[\max\limits_{x\in\bigcup\limits_{j=1}^k w_j\widetilde{\mbox{Poly}}_{f,\mathscr{U}}(m_j,M_j)}\frac{x}{\log x}\right]\log\left(\sum\limits_{j=1}^k w_j\bm{V}^{\ast}\mbox{Poly}_{f,\mathscr{U}}(\bm{A}_j)\bm{V}\right);
\end{eqnarray}
and, the following lower bound for $\sum\limits_{j=1}^k w_j\Phi(f(\bm{A}_j))$:
\begin{eqnarray}\label{eq log x l0 LB: cor: 2.11-14}
\sum\limits_{j=1}^k w_j\Phi(f(\bm{A}_j))&\geq&\left[\min\limits_{x\in\bigcup\limits_{j=1}^k w_j\widetilde{\mbox{Poly}}_{f,\mathscr{L}}(m_j,M_j)}\frac{x}{\log x}\right]\log\left(\sum\limits_{j=1}^k w_j\bm{V}^{\ast}\mbox{Poly}_{f,\mathscr{L}}(\bm{A}_j)\bm{V}\right).
\end{eqnarray}

(II') If we also assume that $\log\left(\sum\limits_{j=1}^k w_j\bm{V}^{\ast}\mbox{Poly}_{f,\mathscr{U}}(\bm{A}_j)\bm{V}\right) < \bm{O}$ and $\log\left(\sum\limits_{j=1}^k w_j\bm{V}^{\ast}\mbox{Poly}_{f,\mathscr{L}}(\bm{A}_j)\bm{V}\right) < \bm{O}$, then, we have the following upper bound for $\sum\limits_{j=1}^k w_j\Phi(f(\bm{A}_j))$:
\begin{eqnarray}\label{eq log x s0 UB: cor: 2.11-14}
\sum\limits_{j=1}^k w_j\Phi(f(\bm{A}_j))&\leq&\left[\min\limits_{x\in\bigcup\limits_{j=1}^k w_j\widetilde{\mbox{Poly}}_{f,\mathscr{L}}(m_j,M_j)}\frac{x}{\log x}\right]\log\left(\sum\limits_{j=1}^k w_j\bm{V}^{\ast}\mbox{Poly}_{f,\mathscr{L}}(\bm{A}_j)\bm{V}\right);
\end{eqnarray}
and, the following lower bound for $\sum\limits_{j=1}^k w_j\Phi(f(\bm{A}_j))$:
\begin{eqnarray}\label{eq log x s0 LB: cor: 2.11-14}
\sum\limits_{j=1}^k w_j\Phi(f(\bm{A}_j))&\geq&\left[\max\limits_{x\in\bigcup\limits_{j=1}^k w_j\widetilde{\mbox{Poly}}_{f,\mathscr{U}}(m_j,M_j)}\frac{x}{\log x}\right]\log\left(\sum\limits_{j=1}^k w_j\bm{V}^{\ast}\mbox{Poly}_{f,\mathscr{U}}(\bm{A}_j)\bm{V}\right).
\end{eqnarray}

(III) Since we have $\exp\left(\sum\limits_{j=1}^k w_j\bm{V}^{\ast}\mbox{Poly}_{f,\mathscr{U}}(\bm{A}_j)\bm{V}\right) > \bm{O}$ and $\exp\left(\sum\limits_{j=1}^k w_j\bm{V}^{\ast}\mbox{Poly}_{f,\mathscr{L}}(\bm{A}_j)\bm{V}\right) > \bm{O}$, then, we have the following upper bound for $\sum\limits_{j=1}^k w_j\Phi(f(\bm{A}_j))$:
\begin{eqnarray}\label{eq exp x UB: cor: 2.11-14}
\sum\limits_{j=1}^k w_j\Phi(f(\bm{A}_j))&\leq&\left[\max\limits_{x\in\bigcup\limits_{j=1}^k w_j\widetilde{\mbox{Poly}}_{f,\mathscr{U}}(m_j,M_j)}\frac{x}{\exp x}\right]\exp\left(\sum\limits_{j=1}^k w_j\bm{V}^{\ast}\mbox{Poly}_{f,\mathscr{U}}(\bm{A}_j)\bm{V}\right);
\end{eqnarray}
and, the following lower bound for $\sum\limits_{j=1}^k w_j\Phi(f(\bm{A}_j))$:
\begin{eqnarray}\label{eq exp x LB: cor: 2.11-14}
\sum\limits_{j=1}^k w_j\Phi(f(\bm{A}_j))&\geq&\left[\min\limits_{x\in\bigcup\limits_{j=1}^k w_j\widetilde{\mbox{Poly}}_{f,\mathscr{L}}(m_j,M_j)}\frac{x}{\exp x}\right]\exp\left(\sum\limits_{j=1}^k w_j\bm{V}^{\ast}\mbox{Poly}_{f,\mathscr{L}}(\bm{A}_j)\bm{V}\right).
\end{eqnarray}
\end{corollary}
\textbf{Proof:}
Part (I) of this Corollary is proved by applying Theorem~\ref{thm: 2.9} Part (I) with the function $g$ as $g(x)=x^q$. Part (II) of this Corollary is proved by applying Theorem~\ref{thm: 2.9} Part (I) with the function $g$ as $g(x)=\log(x)$, where $\log(x) > 0$. Part (II') of this Corollary is proved by applying Theorem~\ref{thm: 2.9} Part (II) with the function $g$ as $g(x)=\log(x)$, where $\log(x) < 0$. Finally, Part (III) of this Corollary is proved by applying Theorem~\ref{thm: 2.9} Part (I) with the function $g$ as $g(x)=\exp(x)$.      
$\hfill \Box$

\section{Generalized Converses of Operator Jensen’s Inequalities: Difference Kind}\label{sec: Generalized Converses of Operator Jensen’s Inequalities: Difference Kind}

In this section, we will derive the lower and upper bounds for $\sum\limits_{j=1}^k w_j\Phi(f(\bm{A}_j))$ in terms of difference criteria related to the function $g$.
\begin{theorem}\label{thm: cor 2.15}
Let $\bm{A}_j$ be self-adjoint operator with $\Lambda(\bm{A}_j) \in [m_j, M_j]$ for real scalars $m_j <  M_j$. The mapping $\Phi: \mathbb{B}(\mathfrak{H}) \rightarrow \mathbb{B}(\mathfrak{K})$ is defined by Eq.~\eqref{eq: new phi def}. The index $j$ is in the range of $1,2,\cdots,k$, and we have a probability vector $\bm{w}=[w_1,w_2,\cdots, w_k]$ with the dimension $k$, i.e., $\sum\limits_{j=1}^{k}w_j = 1$. Let $f$ be any real valued continous functions defined on the range $\bigcup\limits_{j=1}^k [m_j, M_j]$, represented by $f \in \mathcal{C}(\bigcup\limits_{j=1}^k [m_j, M_j])$. Besides, we assume that the function $f$ satisfes thsoe conditions provided by Eq.~\eqref{eq1: lma: lower and upper bound for f(A)} and Eq.~\eqref{eq: lower poly formats}. The function $g$ is also a real valued continous function defined on the range $\left(\bigcup\limits_{j=1}^k w_j\widetilde{\mbox{Poly}}_{f,\mathscr{L}}(m_j,M_j)\right)\bigcup\left(\bigcup\limits_{j=1}^k w_j\widetilde{\mbox{Poly}}_{f,\mathscr{U}}(m_j,M_j)\right)$. We also have a real valued function $F(u,v)$ defined as Eq.~\eqref{eq: F u v } with support domain on $U \times V$ such that $f(\bigcup\limits_{j=1}^k [m_j, M_j]) \subset U$, and  $g\left(\left(\bigcup\limits_{j=1}^k w_j\widetilde{\mbox{Poly}}_{f,\mathscr{L}}(m_j,M_j)\right)\bigcup\left(\bigcup\limits_{j=1}^k w_j\widetilde{\mbox{Poly}}_{f,\mathscr{U}}(m_j,M_j)\right)\right) \subset V$. 

Then, we have the following upper bound:
\begin{eqnarray}\label{eq UB: cor 2.15}
\sum\limits_{j=1}^k w_j\Phi(f(\bm{A}_j)) - g\left(\sum\limits_{j=1}^k w_j\bm{V}^{\ast}\mbox{Poly}_{f,\mathscr{U}}(\bm{A}_j)\bm{V}\right) \leq 
\max\limits_{x\in\bigcup\limits_{j=1}^k w_j\widetilde{\mbox{Poly}}_{f,\mathscr{U}}(m_j,M_j)}(x - g(x))\bm{1}_{\mathfrak{K}}.
\end{eqnarray}
Similarly, we also have the following lower bound:
\begin{eqnarray}\label{eq LB: cor 2.15}
\sum\limits_{j=1}^k w_j\Phi(f(\bm{A}_j)) - g\left(\sum\limits_{j=1}^k w_j\bm{V}^{\ast}\mbox{Poly}_{f,\mathscr{L}}(\bm{A}_j)\bm{V}\right) \geq
\min\limits_{x\in\bigcup\limits_{j=1}^k w_j\widetilde{\mbox{Poly}}_{f,\mathscr{L}}(m_j,M_j)}(x - g(x))\bm{1}_{\mathfrak{K}}.
\end{eqnarray}
\end{theorem}
\textbf{Proof:}
The upper bound of this theorem is proved by setting $\alpha=1$ in Eq.~\eqref{eq UB: thm: 2.4} from Theorem~\ref{thm: 2.4} and rearrangement of the term $ g\left(\sum\limits_{j=1}^k w_j\bm{V}^{\ast}\mbox{Poly}_{f,\mathscr{U}}(\bm{A}_j)\bm{V}\right)$ to obtain Eq.~\eqref{eq UB: cor 2.15}. 

Similarly, the lower bound of this theorem is proved by setting $\alpha=1$ in Eq.~\eqref{eq LB: thm: 2.4} from Theorem~\ref{thm: 2.4} and rearrangement of the term $ g\left(\sum\limits_{j=1}^k w_j\bm{V}^{\ast}\mbox{Poly}_{f,\mathscr{L}}(\bm{A}_j)\bm{V}\right)$ to obtain Eq.~\eqref{eq LB: cor 2.15}. 
$\hfill \Box$

Next Corollary~\ref{cor: 2.17 ext} is obtained by applying Theorem~\ref{thm: cor 2.15} to special types of function $g$. 
\begin{corollary}\label{cor: 2.17 ext}
Let $\bm{A}_j$ be self-adjoint operator with $\Lambda(\bm{A}_j) \in [m_j, M_j]$ for real scalars $m_j <  M_j$. The mapping $\Phi: \mathbb{B}(\mathfrak{H}) \rightarrow \mathbb{B}(\mathfrak{K})$ is defined by Eq.~\eqref{eq: new phi def}. The index $j$ is in the range of $1,2,\cdots,k$, and we have a probability vector $\bm{w}=[w_1,w_2,\cdots, w_k]$ with the dimension $k$, i.e., $\sum\limits_{j=1}^{k}w_j = 1$. Let $f$ be any real valued continous functions defined on the range $\bigcup\limits_{j=1}^k [m_j, M_j]$, represented by $f \in \mathcal{C}(\bigcup\limits_{j=1}^k [m_j, M_j])$. Besides, we assume that the function $f$ satisfes thsoe conditions provided by Eq.~\eqref{eq1: lma: lower and upper bound for f(A)} and Eq.~\eqref{eq: lower poly formats}. 

(I) If we also assume that $g(x)=x^q$, where $q \in \mathbb{R}$, then, we have the following upper bound for $\sum\limits_{j=1}^k w_j\Phi(f(\bm{A}_j))$:
\begin{eqnarray}\label{eq x q UB: cor: 2.17 ext}
\sum\limits_{j=1}^k w_j\Phi(f(\bm{A}_j)) - \left(\sum\limits_{j=1}^k w_j\bm{V}^{\ast}\mbox{Poly}_{f,\mathscr{U}}(\bm{A}_j)\bm{V}\right)^q \leq 
\max\limits_{x\in\bigcup\limits_{j=1}^k w_j\widetilde{\mbox{Poly}}_{f,\mathscr{U}}(m_j,M_j)}(x - x^q)\bm{1}_{\mathfrak{K}};
\end{eqnarray}
and, the following lower bound for $\sum\limits_{j=1}^k w_j\Phi(f(\bm{A}_j))$:
\begin{eqnarray}\label{eq x q LB: cor: 2.17 ext}
\sum\limits_{j=1}^k w_j\Phi(f(\bm{A}_j)) - \left(\sum\limits_{j=1}^k w_j\bm{V}^{\ast}\mbox{Poly}_{f,\mathscr{L}}(\bm{A}_j)\bm{V}\right)^q \geq
\min\limits_{x\in\bigcup\limits_{j=1}^k w_j\widetilde{\mbox{Poly}}_{f,\mathscr{L}}(m_j,M_j)}(x - x^q)\bm{1}_{\mathfrak{K}}.
\end{eqnarray}

(II) If we also assume that $g(x)=\log(x)$, then, we have the following upper bound for $\sum\limits_{j=1}^k w_j\Phi(f(\bm{A}_j))$:
\begin{eqnarray}\label{eq log x UB: cor: 2.17 ext}
\sum\limits_{j=1}^k w_j\Phi(f(\bm{A}_j)) - \log\left(\sum\limits_{j=1}^k w_j\bm{V}^{\ast}\mbox{Poly}_{f,\mathscr{U}}(\bm{A}_j)\bm{V}\right) \leq 
\max\limits_{x\in\bigcup\limits_{j=1}^k w_j\widetilde{\mbox{Poly}}_{f,\mathscr{U}}(m_j,M_j)}(x - \log(x))\bm{1}_{\mathfrak{K}};
\end{eqnarray}
and, the following lower bound for $\sum\limits_{j=1}^k w_j\Phi(f(\bm{A}_j))$:
\begin{eqnarray}\label{eq log x LB: cor: 2.17 ext}
\sum\limits_{j=1}^k w_j\Phi(f(\bm{A}_j)) - \log\left(\sum\limits_{j=1}^k w_j\bm{V}^{\ast}\mbox{Poly}_{f,\mathscr{L}}(\bm{A}_j)\bm{V}\right) \geq
\min\limits_{x\in\bigcup\limits_{j=1}^k w_j\widetilde{\mbox{Poly}}_{f,\mathscr{L}}(m_j,M_j)}(x - \log(x))\bm{1}_{\mathfrak{K}}.
\end{eqnarray}

(III) If we also assume that $g(x)=\exp(x)$, then, we have the following upper bound for $\sum\limits_{j=1}^k w_j\Phi(f(\bm{A}_j))$:
\begin{eqnarray}\label{eq exp x UB: cor: 2.17 ext}
\sum\limits_{j=1}^k w_j\Phi(f(\bm{A}_j)) - \exp\left(\sum\limits_{j=1}^k w_j\bm{V}^{\ast}\mbox{Poly}_{f,\mathscr{U}}(\bm{A}_j)\bm{V}\right) \leq 
\max\limits_{x\in\bigcup\limits_{j=1}^k w_j\widetilde{\mbox{Poly}}_{f,\mathscr{U}}(m_j,M_j)}(x - \exp(x))\bm{1}_{\mathfrak{K}};
\end{eqnarray}
and, the following lower bound for $\sum\limits_{j=1}^k w_j\Phi(f(\bm{A}_j))$:
\begin{eqnarray}\label{eq exp x LB: cor: 2.17 ext}
\sum\limits_{j=1}^k w_j\Phi(f(\bm{A}_j)) - \exp\left(\sum\limits_{j=1}^k w_j\bm{V}^{\ast}\mbox{Poly}_{f,\mathscr{L}}(\bm{A}_j)\bm{V}\right) \geq
\min\limits_{x\in\bigcup\limits_{j=1}^k w_j\widetilde{\mbox{Poly}}_{f,\mathscr{L}}(m_j,M_j)}(x - \exp(x))\bm{1}_{\mathfrak{K}}.
\end{eqnarray}
\end{corollary}
\textbf{Proof:}
For Part (I), we will use $g(x)=x^q$ in Eq.~\eqref{eq UB: cor 2.15} in Theorem~\ref{thm: cor 2.15} to obtain Eq.~\eqref{eq x q UB: cor: 2.17 ext}. We will also use $g(x)=x^q$ in Eq.~\eqref{eq LB: cor 2.15} in Theorem~\ref{thm: cor 2.15} to obtain Eq.~\eqref{eq x q LB: cor: 2.17 ext}. 

For Part (II), we will use $g(x)=\log(x)$ in Eq.~\eqref{eq UB: cor 2.15} in Theorem~\ref{thm: cor 2.15} to obtain Eq.~\eqref{eq log x UB: cor: 2.17 ext}. We will also use $g(x)=\log(x)$ in Eq.~\eqref{eq LB: cor 2.15} in Theorem~\ref{thm: cor 2.15} to obtain Eq.~\eqref{eq log x LB: cor: 2.17 ext}. 

For Part (III), we will use $g(x)=\exp(x)$ in Eq.~\eqref{eq UB: cor 2.15} in Theorem~\ref{thm: cor 2.15} to obtain Eq.~\eqref{eq exp x UB: cor: 2.17 ext}. We will also use $g(x)=\exp(x)$ in Eq.~\eqref{eq LB: cor 2.15} in Theorem~\ref{thm: cor 2.15} to obtain Eq.~\eqref{eq exp x LB: cor: 2.17 ext}. 
$\hfill \Box$

\section{Hypercomplex Function Approximation}\label{sec: Hypercomplex Function Approximations}

In this section, we will consider hypercomplex function approximation problem in terms of the ratio error discussed in Section~\ref{sec: Ratio Type Approximation} and the difference error in Section~\ref{sec: Difference Type Approximation}.

\subsection{Ratio Type Approximation}\label{sec: Ratio Type Approximation}

Let $\bm{A}$ be a self-adjoint operator with $\Lambda(\bm{A}) \in [m, M]$ for real scalars $m <  M$. The mapping $\Phi: \mathbb{B}(\mathfrak{H}) \rightarrow \mathbb{B}(\mathfrak{K})$ is defined by Eq.~\eqref{eq: new phi def}. Let $f$ be any real valued continuous functions defined on the range $[m, M]$, represented by $f \in \mathcal{C}([m, M])$. The \emph{ratio type approximation} problem is to find the function $g$ and the polynomial function $p_1$ to satisfy the following:
\begin{eqnarray}\label{eq1: UP Ratio Type Approximation}
\frac{\Phi(f(\bm{A}))}{g(\bm{V}^{\ast}p_1(\bm{A})\bm{V})}\leq \alpha_1 \bm{1}_{\mathfrak{K}},
\end{eqnarray}
where $\alpha_1$ is some specified positive real number. Similarly, we also can find the polynomial function $p_2$ to satisfy the following:
\begin{eqnarray}\label{eq1: LO Ratio Type Approximation}
\frac{\Phi(f(\bm{A}))}{g(\bm{V}^{\ast}p_2(\bm{A})\bm{V})}\geq \alpha_2 \bm{1}_{\mathfrak{K}},
\end{eqnarray}
where $\alpha_2$ is some specified positive real number.

From Theorem~\ref{thm: 2.9}, if $g\left(\bm{V}^{\ast}\mbox{Poly}_{f,\mathscr{U}}(\bm{A})\bm{V}\right) > \bm{O}$, Eq.~\eqref{eq1: UP Ratio Type Approximation} can be established by setting 
\begin{eqnarray}\label{eq2: UP Ratio Type Approximation}
\alpha_1 &\geq& \left[\max\limits_{x\in\widetilde{\mbox{Poly}}_{f,\mathscr{U}}(m,M)}xg^{-1}(x)\right]. 
\end{eqnarray}
Similarly,  if $g\left(\bm{V}^{\ast}\mbox{Poly}_{f,\mathscr{L}}(\bm{A})\bm{V}\right) > \bm{O}$, Eq.~\eqref{eq1: LO Ratio Type Approximation} can be established by setting 
\begin{eqnarray}\label{eq2: LO Ratio Type Approximation}
\alpha_2 &\leq& \left[\min\limits_{x\in\widetilde{\mbox{Poly}}_{f,\mathscr{L}}(m,M)}xg^{-1}(x)\right]. 
\end{eqnarray}
On the other hand, from Theorem~\ref{thm: 2.9}, if $g\left(\bm{V}^{\ast}\mbox{Poly}_{f,\mathscr{U}}(\bm{A})\bm{V}\right) < \bm{O}$, Eq.~\eqref{eq1: UP Ratio Type Approximation} can be established by setting 
\begin{eqnarray}\label{eq2-1: UP Ratio Type Approximation}
\alpha_1 &\geq& \left[\min\limits_{x\in\widetilde{\mbox{Poly}}_{f,\mathscr{L}}(m,M)}xg^{-1}(x)\right]. 
\end{eqnarray}
Similarly,  if $g\left(\bm{V}^{\ast}\mbox{Poly}_{f,\mathscr{L}}(\bm{A})\bm{V}\right) < \bm{O}$, Eq.~\eqref{eq1: LO Ratio Type Approximation} can be established by setting 
\begin{eqnarray}\label{eq2-1: LO Ratio Type Approximation}
\alpha_2 &\leq& \left[\max\limits_{x\in\widetilde{\mbox{Poly}}_{f,\mathscr{U}}(m,M)}xg^{-1}(x)\right]. 
\end{eqnarray}

For $k$ terms, let $\bm{A}_j$ be self-adjoint operator with $\Lambda(\bm{A}_j) \in [m_j, M_j]$ for real scalars $m_j <  M_j$. The mapping $\Phi: \mathbb{B}(\mathfrak{H}) \rightarrow \mathbb{B}(\mathfrak{K})$ is defined by Eq.~\eqref{eq: new phi def}. The index $j$ is in the range of $1,2,\cdots,k$, and we have a probability vector $\bm{w}=[w_1,w_2,\cdots, w_k]$ with the dimension $k$, i.e., $\sum\limits_{j=1}^{k}w_j = 1$. Let $f$ be any real valued continuous functions defined on the range $\bigcup\limits_{j=1}^k [m_j, M_j]$, represented by $f \in \mathcal{C}(\bigcup\limits_{j=1}^k [m_j, M_j])$. Besides, we assume that the function $f$ satisfies the conditions provided by Eq.~\eqref{eq1: lma: lower and upper bound for f(A)} and Eq.~\eqref{eq: lower poly formats}. The \emph{ratio type approximation} problem is to find the function $g$ and the polynomial functions $p_{1,1}(x),\cdots,p_{1,k}(x)$ to satisfy the following:
\begin{eqnarray}\label{eq3: UP Ratio Type Approximation}
\frac{\sum\limits_{j=1}^k w_j\Phi(f(\bm{A}_j))}{g\left(\sum\limits_{j=1}^k w_j\bm{V}^{\ast}p_{1,j}(\bm{A}_j)\bm{V}\right)}\leq\alpha_1 \bm{1}_{\mathfrak{K}},
\end{eqnarray}
where $\alpha_1$ is some specfied positive real number. Similarly, we also can find the polynomial function $p_{2,1}(x),\cdots,p_{2,k}(x)$  to satisfy the following:
\begin{eqnarray}\label{eq3: LO Ratio Type Approximation}
\frac{\sum\limits_{j=1}^k w_j\Phi(f(\bm{A}_j))}{g\left(\sum\limits_{j=1}^k w_j\bm{V}^{\ast}p_{2,j}(\bm{A}_j)\bm{V}\right)}\geq \alpha_2 \bm{1}_{\mathfrak{K}},
\end{eqnarray}
where $\alpha_2$ is some specfied positive real number. 

From Theorem~\ref{thm: 2.9}, if $g\left(\bm{V}^{\ast}\mbox{Poly}_{f,\mathscr{U}}(\bm{A})\bm{V}\right) > \bm{O}$, Eq.~\eqref{eq3: UP Ratio Type Approximation} can be established by setting 
\begin{eqnarray}\label{eq4: UP Ratio Type Approximation}
\alpha_1 &\geq& \left[\max\limits_{x\in\bigcup\limits_{j=1}^k w_j\widetilde{\mbox{Poly}}_{f,\mathscr{U}}(m_j,M_j)}xg^{-1}(x)\right]. 
\end{eqnarray}
Similarly,  if $g\left(\bm{V}^{\ast}\mbox{Poly}_{f,\mathscr{L}}(\bm{A})\bm{V}\right) > \bm{O}$, Eq.~\eqref{eq3: LO Ratio Type Approximation} can be established by setting 
\begin{eqnarray}\label{eq4: LO Ratio Type Approximation}
\alpha_2 &\leq& \left[\min\limits_{x\in\bigcup\limits_{j=1}^k w_j\widetilde{\mbox{Poly}}_{f,\mathscr{L}}(m_j,M_j)}xg^{-1}(x)\right]. 
\end{eqnarray}
On the other hand, from Theorem~\ref{thm: 2.9}, if $g\left(\bm{V}^{\ast}\mbox{Poly}_{f,\mathscr{U}}(\bm{A})\bm{V}\right) < \bm{O}$, Eq.~\eqref{eq3: UP Ratio Type Approximation} can be established by setting 
\begin{eqnarray}\label{eq4-1: UP Ratio Type Approximation}
\alpha_1 &\geq& \left[\min\limits_{x\in\bigcup\limits_{j=1}^k w_j\widetilde{\mbox{Poly}}_{f,\mathscr{L}}(m_j,M_j)}xg^{-1}(x)\right]. 
\end{eqnarray}
Similarly,  if $g\left(\bm{V}^{\ast}\mbox{Poly}_{f,\mathscr{L}}(\bm{A})\bm{V}\right) < \bm{O}$, Eq.~\eqref{eq3: LO Ratio Type Approximation} can be established by setting 
\begin{eqnarray}\label{eq4-1: LO Ratio Type Approximation}
\alpha_2 &\leq& \left[\max\limits_{x\in\bigcup\limits_{j=1}^k w_j\widetilde{\mbox{Poly}}_{f,\mathscr{U}}(m_j,M_j)}xg^{-1}(x)\right]. 
\end{eqnarray}

\begin{remark}\label{rmk: ratio approx conjecture}
Given requirements of Eq.~\eqref{eq1: UP Ratio Type Approximation}, Eq.~\eqref{eq1: LO Ratio Type Approximation}, Eq.~\eqref{eq3: UP Ratio Type Approximation}, and Eq.~\eqref{eq3: LO Ratio Type Approximation}, we conjecture the existence of function $g$ and polynomials $p_1(x), p_2(x), p_{1,1}(x),\cdots,p_{1,k}(x)$, and $p_{2,1}(x),\cdots,p_{2,k}(x)$ for any given $\alpha_1$ or $\alpha_2$. If existence, how to find these functions?
\end{remark}

The following Example~\ref{exp: thm 2.9} is provided to evaluate the upper ratio bound given by Eq.~\eqref{eq3: UP Ratio Type Approximation} and evaluate the lower ratio bound given by Eq.~\eqref{eq3: LO Ratio Type Approximation}. 
\begin{example}\label{exp: thm 2.9}
Let $\bm{A}_j$ be self-adjoint operator with $\Lambda(\bm{A}_j) \in [m, M]$ for real scalars $m <  M$. The mapping $\Phi: \mathbb{B}(\mathfrak{H}) \rightarrow \mathbb{B}(\mathfrak{K})$ is defined as $\Phi(\bm{X})=\bm{V}^{\ast}(\bm{X})\bm{X}$. The index $j$ is in the range of $1,2,\cdots,k$, and we have a probability vector $\bm{w}=[w_1,w_2,\cdots, w_k]$ with the dimension $k$, i.e., $\sum\limits_{j=1}^{k}w_j = 1$. Let $f$ be a convex and differentiable function such that $ax + b'  \leq f(x) \leq ax + b$ for $x \in [m,M]$, where 
\begin{eqnarray}\label{eq1: exp: thm 2.9}
a&=&\frac{f(M) - f(m)}{M-m},\nonumber \\ 
b&=&\frac{Mf(m) - mf(M)}{M-m},\nonumber \\ 
b'&=&f(x_0)-\frac{f(M) - f(m)}{M-m}x_0,
\end{eqnarray}
where $f'(x_0) = \frac{f(M) - f(m)}{M-m}$.

If we require $g(x)>0$ for all $x \in [m,M]$, from Theorem~\ref{thm: 2.9}, we have
\begin{eqnarray}\label{eq2: UP exp: thm 2.9}
\sum\limits_{j=1}^k w_j \Phi(f(\bm{A}_j)) \leq \max\limits_{m<x<M}\frac{ax+b}{g(x)}g\left(\sum\limits_{j=1}^k w_j \Phi(\bm{A}_j)\right),
\end{eqnarray}
and 
\begin{eqnarray}\label{eq2: LO exp: thm 2.9}
\sum\limits_{j=1}^k w_j \Phi(f(\bm{A}_j)) \geq \min\limits_{m<x<M}\frac{ax+b'}{g(x)}g\left(\sum\limits_{j=1}^k w_j \Phi(\bm{A}_j)\right).
\end{eqnarray}

Let us consider several special cases of the function $g(x)$. If $g(x)=x^q$, where $q \in \mathbb{R}$, and $m>0$, we have
\begin{eqnarray}\label{eq2-1: exp: thm 2.9}
\alpha_1 &\geq& \max\limits_{m<x<M}\frac{ax+b}{x^q},\nonumber \\
\alpha_2 &\leq& \min\limits_{m<x<M}\frac{ax+b'}{x^q}.
\end{eqnarray}
If $g(x)=\log(x)$ and $m>1$, we have
\begin{eqnarray}\label{eq2-2: exp: thm 2.9}
\alpha_1 &\geq& \max\limits_{m<x<M}\frac{ax+b}{\log(x)},\nonumber \\
\alpha_2 &\leq& \min\limits_{m<x<M}\frac{ax+b'}{\log(x)}.
\end{eqnarray}
If $g(x)=\exp(x)$, we have
\begin{eqnarray}\label{eq2-2: exp: thm 2.9}
\alpha_1 &\geq& \max\limits_{m<x<M}\frac{ax+b}{\exp(x)},\nonumber \\
\alpha_2 &\leq& \min\limits_{m<x<M}\frac{ax+b'}{\exp(x)}.
\end{eqnarray}

If we require $g(x)<0$ for all $x \in [m,M]$, from Theorem~\ref{thm: 2.9}, we also have
\begin{eqnarray}\label{eq3: UP exp: thm 2.9}
\sum\limits_{j=1}^k w_j \Phi(f(\bm{A}_j)) \leq \min\limits_{m<x<M}\frac{ax+b}{g(x)},
\end{eqnarray}
and 
\begin{eqnarray}\label{eq3: LO exp: thm 2.9}
\sum\limits_{j=1}^k w_j \Phi(f(\bm{A}_j)) \geq \max\limits_{m<x<M}\frac{ax+b'}{g(x)},
\end{eqnarray}

Let us consider several special cases of the function $g(x)$. If $g(x)=-x^q$, where $q \in \mathbb{R}$, and $m>0$, we have
\begin{eqnarray}\label{eq3-1: exp: thm 2.9}
\alpha_1 &\geq& \min\limits_{m<x<M}-\frac{ax+b}{x^q},\nonumber \\
\alpha_2 &\leq& \max\limits_{m<x<M}-\frac{ax+b'}{x^q}.
\end{eqnarray}
If $g(x)=\log(x)$ and $0<m<M<1$, we have
\begin{eqnarray}\label{eq3-2: exp: thm 2.9}
\alpha_1 &\geq& \min\limits_{m<x<M}\frac{ax+b}{\log(x)},\nonumber \\
\alpha_2 &\leq& \max\limits_{m<x<M}\frac{ax+b'}{\log(x)}.
\end{eqnarray}
\end{example}

\subsection{Difference Type Approximation}\label{sec: Difference Type Approximation}

Let $\bm{A}$ be a self-adjoint operator with $\Lambda(\bm{A}) \in [m, M]$ for real scalars $m <  M$. The mapping $\Phi: \mathbb{B}(\mathfrak{H}) \rightarrow \mathbb{B}(\mathfrak{K})$ is defined by Eq.~\eqref{eq: new phi def}. Let $f$ be any real valued continuous functions defined on the range $[m, M]$, represented by $f \in \mathcal{C}([m, M])$. The \emph{difference type approximation} problem is to find the function $g$ and the polynomial function $p_1(x)$ to satisfy the following:
\begin{eqnarray}\label{eq1: UP Difference Type Approximation}
\Phi(f(\bm{A}))- g(\bm{V}^{\ast}p_1(\bm{A})\bm{V})\leq \beta_1 \bm{1}_{\mathfrak{K}},
\end{eqnarray}
where $\beta_1$ is some specified real number. Similarly, we also can find the polynomial function $p_2(x)$ to satisfy the following:
\begin{eqnarray}\label{eq1: LO Difference Type Approximation}
\Phi(f(\bm{A}))-g(\bm{V}^{\ast}p_2(\bm{A})\bm{V})\geq \beta_2 \bm{1}_{\mathfrak{K}},
\end{eqnarray}
where $\beta_2$ is some specified real number.

From Theorem~\ref{thm: cor 2.15}, Eq.~\eqref{eq1: UP Difference Type Approximation} can be established by setting 
\begin{eqnarray}\label{eq2: UP Difference Type Approximation}
\beta_1 &\geq& \left[\max\limits_{x\in\widetilde{\mbox{Poly}}_{f,\mathscr{U}}(m,M)}x-g(x)\right]. 
\end{eqnarray}
Similarly, Eq.~\eqref{eq1: LO Difference Type Approximation} can be established by setting 
\begin{eqnarray}\label{eq2: LO Difference Type Approximation}
\beta_2 &\leq& \left[\min\limits_{x\in\widetilde{\mbox{Poly}}_{f,\mathscr{L}}(m,M)}x - g(x)\right]. 
\end{eqnarray}

For $k$ terms, let $\bm{A}_j$ be self-adjoint operator with $\Lambda(\bm{A}_j) \in [m_j, M_j]$ for real scalars $m_j <  M_j$. The mapping $\Phi: \mathbb{B}(\mathfrak{H}) \rightarrow \mathbb{B}(\mathfrak{K})$ is defined by Eq.~\eqref{eq: new phi def}. The index $j$ is in the range of $1,2,\cdots,k$, and we have a probability vector $\bm{w}=[w_1,w_2,\cdots, w_k]$ with the dimension $k$, i.e., $\sum\limits_{j=1}^{k}w_j = 1$. Let $f$ be any real valued continuous functions defined on the range $\bigcup\limits_{j=1}^k [m_j, M_j]$, represented by $f \in \mathcal{C}(\bigcup\limits_{j=1}^k [m_j, M_j])$. Besides, we assume that the function $f$ satisfies the conditions provided by Eq.~\eqref{eq1: lma: lower and upper bound for f(A)} and Eq.~\eqref{eq: lower poly formats}. The \emph{difference type approximation} problem is to find the function $g$ and the polynomial functions $p_{1,1}(x),\cdots,p_{1,k}(x)$ to satisfy the following:
\begin{eqnarray}\label{eq3: UP Difference Type Approximation}
\sum\limits_{j=1}^k w_j\Phi(f(\bm{A}_j))-g\left(\sum\limits_{j=1}^k w_j\bm{V}^{\ast}p_{1,j}(\bm{A}_j)\bm{V}\right)\leq\beta_1 \bm{1}_{\mathfrak{K}},
\end{eqnarray}
where $\beta_1$ is some specfied real number. Similarly, we also can find the polynomial function $p_{2,1}(x),\cdots,p_{2,k}(x)$  to satisfy the following:
\begin{eqnarray}\label{eq3: LO Difference Type Approximation}
\sum\limits_{j=1}^k w_j\Phi(f(\bm{A}_j))-g\left(\sum\limits_{j=1}^k w_j\bm{V}^{\ast}p_{2,j}(\bm{A}_j)\bm{V}\right)\geq\beta_2 \bm{1}_{\mathfrak{K}},
\end{eqnarray}
where $\beta_2$ is some specfied real number. 

The following Example~\ref{exp: cor 2.15} is provided to evaluate the upper difference bound given by Eq.~\eqref{eq3: UP Difference Type Approximation} and evaluate the lower difference bound given by Eq.~\eqref{eq3: LO Difference Type Approximation}. 
\begin{example}\label{exp: cor 2.15}
Let $\bm{A}_j$ be self-adjoint operator with $\Lambda(\bm{A}_j) \in [m, M]$ for real scalars $m <  M$. The mapping $\Phi: \mathbb{B}(\mathfrak{H}) \rightarrow \mathbb{B}(\mathfrak{K})$ is defined as $\Phi(\bm{X})=\bm{V}^{\ast}(\bm{X})\bm{X}$. The index $j$ is in the range of $1,2,\cdots,k$, and we have a probability vector $\bm{w}=[w_1,w_2,\cdots, w_k]$ with the dimension $k$, i.e., $\sum\limits_{j=1}^{k}w_j = 1$. Let $f$ be a convex and differentiable function such that $ax + b'  \leq f(x) \leq ax + b$ for $x \in [m,M]$, where 
\begin{eqnarray}\label{eq1: exp: cor 2.15}
a&=&\frac{f(M) - f(m)}{M-m},\nonumber \\ 
b&=&\frac{Mf(m) - mf(M)}{M-m},\nonumber \\ 
b'&=&f(x_0)-\frac{f(M) - f(m)}{M-m}x_0,
\end{eqnarray}
where $f'(x_0) = \frac{f(M) - f(m)}{M-m}$.

From Theorem~\ref{thm: cor 2.15}, we have
\begin{eqnarray}\label{eq2: UP exp: cor 2.15}
\sum\limits_{j=1}^k w_j \Phi(f(\bm{A}_j)) \leq \max\limits_{m<x<M}(ax+b-g(x)),
\end{eqnarray}
and 
\begin{eqnarray}\label{eq2: LO exp: cor 2.15}
\sum\limits_{j=1}^k w_j \Phi(f(\bm{A}_j)) \geq \min\limits_{m<x<M}(ax+b'-g(x)),
\end{eqnarray}

Let us consider several special cases of the function $g(x)$. If $g(x)=x^q$, where $q \in \mathbb{R}$ and $m>0$, we have
\begin{eqnarray}\label{eq2-1: exp: cor 2.15}
\beta_1 &\geq& \max\limits_{m<x<M}(ax+b-x^q),\nonumber \\
\beta_2 &\leq& \min\limits_{m<x<M}(ax+b'-x^q).
\end{eqnarray}
If $g(x)=\log(x)$ and $m>0$, we have
\begin{eqnarray}\label{eq2-2: exp: cor 2.15}
\beta_1 &\geq& \max\limits_{m<x<M}(ax+b-\log(x)),\nonumber \\
\beta_2 &\leq& \min\limits_{m<x<M}(ax+b'-\log(x)).
\end{eqnarray}
If $g(x)=\exp(x)$, we have
\begin{eqnarray}\label{eq2-2: exp: cor 2.15}
\beta_1 &\geq& \max\limits_{m<x<M}(ax+b-\exp(x)),\nonumber \\
\beta_2 &\leq& \min\limits_{m<x<M}(ax+b'-\exp(x)).
\end{eqnarray}
\end{example}

\begin{remark}\label{rmk: difference approx conjecture}
Given requirements of Eq.~\eqref{eq1: UP Difference Type Approximation}, Eq.~\eqref{eq1: LO Difference Type Approximation}, Eq.~\eqref{eq3: UP Difference Type Approximation}, and Eq.~\eqref{eq3: LO Difference Type Approximation}, we conjecture the existence of function $g$ and polynomials $p_1(x), p_2(x), p_{1,1}(x),\cdots,p_{1,k}(x)$, and $p_{2,1}(x),\cdots,p_{2,k}(x)$ for any given $\beta_1$ or $\beta_2$. If existence, how to find these functions?
\end{remark}

\section{Bounds Algebra}\label{sec: Bounds Algebra}

In this section, we will apply results from Section~\ref{sec: Generalized Converses of Operator Jensen’s Inequalities: Ratio Kind} to build bounds algebra for addition and multiplication of hypercomplex functions 
$\sum\limits_{j=1}^k w_j\Phi(f(\bm{A}_j))$ and $\sum\limits_{j=1}^k w_j\Phi(h(\bm{A}_j))$, where functions $f$ and $h$ share some common properties. 

\subsection{Hypercomplex Function Bounds Algebra}\label{sec: Function of Operator Bounds Algebra}

Let $\bm{A}_j$ be self-adjoint operator with $\Lambda(\bm{A}_j) \in [m, M]$ for real scalars $m <  M$. The mapping $\Phi: \mathbb{B}(\mathfrak{H}) \rightarrow \mathbb{B}(\mathfrak{K})$ is defined as $\Phi(\bm{X})=\bm{V}^{\ast}(\bm{X})\bm{X}$. The index $j$ is in the range of $1,2,\cdots,k$, and we have a probability vector $\bm{w}=[w_1,w_2,\cdots, w_k]$ with the dimension $k$, i.e., $\sum\limits_{j=1}^{k}w_j = 1$. Let $f, h$ be two convex and differentiable functions defined in $[m,M]$ such that 
\begin{eqnarray}\label{eq1: sec: Function of Operator Bounds Algebra}
ax+b'\leq f(x)  \leq ax+b, \nonumber \\
cx+d'\leq h(x) \leq cx+d,
\end{eqnarray}
where $x \in [m,M]$. If we require $g(x)>0$ for all $x \in [m,M]$, from Theorem~\ref{thm: 2.9}, we have
\begin{eqnarray}\label{eq2f: UP sec: Function of Operator Bounds Algebra}
\sum\limits_{j=1}^k w_j \Phi(f(\bm{A}_j)) \leq \underbrace{\max\limits_{m<x<M}\frac{ax+b}{g(x)}}_{:=\alpha_{f,\mathscr{U}}}g\left(\sum\limits_{j=1}^k w_j \Phi(\bm{A}_j)\right),
\end{eqnarray}
and 
\begin{eqnarray}\label{eq2f: LO sec: Function of Operator Bounds Algebra}
\sum\limits_{j=1}^k w_j \Phi(f(\bm{A}_j)) \geq \underbrace{\min\limits_{m<x<M}\frac{ax+b'}{g(x)}}_{:=\alpha_{f,\mathscr{L}}}g\left(\sum\limits_{j=1}^k w_j \Phi(\bm{A}_j)\right).
\end{eqnarray}
Similarly,  from Theorem~\ref{thm: 2.9}, we also have
\begin{eqnarray}\label{eq2h: UP sec: Function of Operator Bounds Algebra}
\sum\limits_{j=1}^k w_j \Phi(h(\bm{A}_j)) \leq \underbrace{\max\limits_{m<x<M}\frac{cx+d}{g(x)}}_{:=\alpha_{h,\mathscr{U}}}g\left(\sum\limits_{j=1}^k w_j \Phi(\bm{A}_j)\right),
\end{eqnarray}
and 
\begin{eqnarray}\label{eq2h: LO sec: Function of Operator Bounds Algebra}
\sum\limits_{j=1}^k w_j \Phi(h(\bm{A}_j)) \geq \underbrace{\min\limits_{m<x<M}\frac{cx+d'}{g(x)}}_{:=\alpha_{h,\mathscr{L}}}g\left(\sum\limits_{j=1}^k w_j \Phi(\bm{A}_j)\right).
\end{eqnarray}

Consider the addition between $\sum\limits_{j=1}^k w_j \Phi(f(\bm{A}_j))$ and $\sum\limits_{j=1}^k w_j \Phi(h(\bm{A}_j))$, from Eq.~\eqref{eq2f: UP sec: Function of Operator Bounds Algebra} and Eq.~\eqref{eq2h: UP sec: Function of Operator Bounds Algebra}, then, we have
\begin{eqnarray}\label{eq2f h: UP sec: Function of Operator Bounds Algebra} 
\sum\limits_{j=1}^k w_j \Phi(f(\bm{A}_j))+\sum\limits_{j=1}^k w_j \Phi(h(\bm{A}_j)) \leq (\alpha_{f,\mathscr{U}} + \alpha_{h,\mathscr{U}})g\left(\sum\limits_{j=1}^k w_j \Phi(\bm{A}_j)\right),
\end{eqnarray}
and, from Eq.~\eqref{eq2f: LO sec: Function of Operator Bounds Algebra} and Eq.~\eqref{eq2h: LO sec: Function of Operator Bounds Algebra}, we also have
\begin{eqnarray}\label{eq2f h: LO sec: Function of Operator Bounds Algebra} 
\sum\limits_{j=1}^k w_j \Phi(f(\bm{A}_j))+\sum\limits_{j=1}^k w_j \Phi(h(\bm{A}_j)) \geq (\alpha_{f,\mathscr{L}} + \alpha_{h,\mathscr{L}})g\left(\sum\limits_{j=1}^k w_j \Phi(\bm{A}_j)\right).
\end{eqnarray}

If any two functions $f,h$ satisfying Eq.~\eqref{eq1: sec: Function of Operator Bounds Algebra}, from Eq.~\eqref{eq2f h: UP sec: Function of Operator Bounds Algebra} and Eq.~\eqref{eq2f h: LO sec: Function of Operator Bounds Algebra}, we can bound the addition between $\sum\limits_{j=1}^k w_j \Phi(f(\bm{A}_j))$ and $\sum\limits_{j=1}^k w_j \Phi(h(\bm{A}_j))$ by coefficients of $g\left(\sum\limits_{j=1}^k w_j \Phi(\bm{A}_j)\right)$ with the term $g\left(\sum\limits_{j=1}^k w_j \Phi(\bm{A}_j)\right)$, which is independent of the functions $f$ and $h$. Therefore, the bounding coefficients of $g\left(\sum\limits_{j=1}^k w_j \Phi(\bm{A}_j)\right)$ form an \emph{algebraic system of interval numbers}, which is an \emph{abelian monoid} with respect to the opration \emph{addition} accoding to Theorem 2.14~\cite{dawood2011theories}. 

Consider the multiplication between $\sum\limits_{j=1}^k w_j \Phi(f(\bm{A}_j))$ and $\sum\limits_{j=1}^k w_j \Phi(h(\bm{A}_j))$, from Eq.~\eqref{eq2f: UP sec: Function of Operator Bounds Algebra} and Eq.~\eqref{eq2h: UP sec: Function of Operator Bounds Algebra} with assumptions of positive $\alpha_{f,\mathscr{U}}, \alpha_{h,\mathscr{U}}, \alpha_{f,\mathscr{L}}$ and $\alpha_{h,\mathscr{L}}$ and $g\left(\sum\limits_{j=1}^k w_j \Phi(\bm{A}_j)\right) > \bm{O}$ , then, we have
\begin{eqnarray}\label{eq3f h: UP sec: Function of Operator Bounds Algebra} 
\left(\sum\limits_{j=1}^k w_j \Phi(f(\bm{A}_j))\right)\times\left(\sum\limits_{j=1}^k w_j \Phi(h(\bm{A}_j))\right)\leq (\alpha_{f,\mathscr{U}} \times \alpha_{h,\mathscr{U}})\left(g\left(\sum\limits_{j=1}^k w_j \Phi(\bm{A}_j)\right)\right)^2,
\end{eqnarray}
and, from Eq.~\eqref{eq2f: LO sec: Function of Operator Bounds Algebra} and Eq.~\eqref{eq2h: LO sec: Function of Operator Bounds Algebra}, we also have
\begin{eqnarray}\label{eq3f h: LO sec: Function of Operator Bounds Algebra} 
\left(\sum\limits_{j=1}^k w_j \Phi(f(\bm{A}_j))\right)\times\left(\sum\limits_{j=1}^k w_j \Phi(h(\bm{A}_j))\right)\geq (\alpha_{f,\mathscr{L}} \times \alpha_{h,\mathscr{L}})\left(g\left(\sum\limits_{j=1}^k w_j \Phi(\bm{A}_j)\right)\right)^2.
\end{eqnarray}

Analogly, if any two functions $f,h$ satisfying Eq.~\eqref{eq1: sec: Function of Operator Bounds Algebra}, from Eq.~\eqref{eq3f h: UP sec: Function of Operator Bounds Algebra} and Eq.~\eqref{eq3f h: LO sec: Function of Operator Bounds Algebra}, we can bound the multiplication between $\sum\limits_{j=1}^k w_j \Phi(f(\bm{A}_j))$ and $\sum\limits_{j=1}^k w_j \Phi(h(\bm{A}_j))$ by coefficients of $\left(g\left(\sum\limits_{j=1}^k w_j \Phi(\bm{A}_j)\right)\right)^2$ with the term $\left(g\left(\sum\limits_{j=1}^k w_j \Phi(\bm{A}_j)\right)\right)^2$, which is independent of the functions $f$ and $h$. Therefore, the bounding coefficients of $\left(g\left(\sum\limits_{j=1}^k w_j \Phi(\bm{A}_j)\right)\right)^2$ form an \emph{algebraic system of interval numbers}, which is an \emph{abelian monoid} with respect to the opration \emph{multiplication} accoding to Theorem 2.14~\cite{dawood2011theories}. 

\begin{remark}\label{rmk: norm bounds}
The norms of the addition (or multiplication) between $\sum\limits_{j=1}^k w_j \Phi(f(\bm{A}_j))$ and $\sum\limits_{j=1}^k w_j \Phi(h(\bm{A}_j))$ can also be bounded by coefficients of $g\left(\sum\limits_{j=1}^k w_j \Phi(\bm{A}_j)\right)$ (or $\left(g\left(\sum\limits_{j=1}^k w_j \Phi(\bm{A}_j)\right)\right)^2$) with abelian monoid algebraic structure and norms of $g\left(\sum\limits_{j=1}^k w_j \Phi(\bm{A}_j)\right)$ (or $\left(g\left(\sum\limits_{j=1}^k w_j \Phi(\bm{A}_j)\right)\right)^2$).
\end{remark}

\subsection{Random Tensor Tail Bounds Algebra}\label{sec: Random Tensor Tail Bounds Algebra}

The purpose of this section is to show a method to obtain tail bound for the addition (or multiplication) between two random tensors $\sum\limits_{j=1}^k w_j \Phi(f(\bm{A}_j))$ and $\sum\limits_{j=1}^k w_j \Phi(h(\bm{A}_j))$ via Theorem~\ref{thm: 2.9}.

Let $\bm{A}_j$ be random Hermitian tensor with $\Lambda(\bm{A}_j) \in [m, M]$ for real scalars $m <  M$. The mapping $\Phi: \mathbb{B}(\mathfrak{H}) \rightarrow \mathbb{B}(\mathfrak{K})$ is defined as $\Phi(\bm{X})=\bm{X}$. The index $j$ is in the range of $1,2,\cdots,k$, and we have a probability vector $\bm{w}=[w_1,w_2,\cdots, w_k]$ with the dimension $k$, i.e., $\sum\limits_{j=1}^{k}w_j = 1$. Let $f, h$ be two convex and differentiable functions defined in $[m,M]$ such that 
\begin{eqnarray}\label{eq1: sec: Random Tensor Tail Bounds Algebra}
f(x)  \leq ax+b, \nonumber \\
h(x) \leq cx+d,
\end{eqnarray}
where $x \in [m,M]$. If we require $g(x)>0$ for all $x \in [m,M]$ and any positive number $\theta$, from Theorem~\ref{thm: 2.9}, we have
\begin{eqnarray}\label{eq2f: UP sec: Random Tensor Tail Bounds Algebra}
\sum\limits_{j=1}^k w_j f(\bm{A}_j) \leq \underbrace{\max\limits_{m<x<M}\frac{ax+b}{g(x)}}_{:=\alpha_{f,\mathscr{U}}}g\left(\sum\limits_{j=1}^k w_j \bm{A}_j\right).
\end{eqnarray}
Similarly,  from Theorem~\ref{thm: 2.9}, we also have
\begin{eqnarray}\label{eq2h: UP sec: Random Tensor Tail Bounds Algebra}
\sum\limits_{j=1}^k w_j h(\bm{A}_j) \leq \underbrace{\min\limits_{m<x<M}\frac{cx+d}{g(x)}}_{:=\alpha_{h,\mathscr{U}}}g\left(\sum\limits_{j=1}^k w_j \bm{A}_j\right).
\end{eqnarray}

Consider tail bounds for the addition between $\sum\limits_{j=1}^k w_j f(\bm{A}_j)$ and $\sum\limits_{j=1}^k w_j h(\bm{A}_j)$, from Eq.~\eqref{eq2f: UP sec: Random Tensor Tail Bounds Algebra} and Eq.~\eqref{eq2h: UP sec: Random Tensor Tail Bounds Algebra}, then, we have
\begin{eqnarray}\label{eq2f h: UP sec: Function of Operator Bounds Algebra}
\mathrm{Pr}\left(\left\Vert\left(\sum\limits_{j=1}^k w_j f(\bm{A}_j) \right)+\left(\sum\limits_{j=1}^k w_j h(\bm{A}_j)\right)\right\Vert_{\ell}\geq\theta\right)\nonumber \\
\leq \mathrm{Pr}\left( \left\Vert(\alpha_{f,\mathscr{U}} + \alpha_{h,\mathscr{U}})g\left(\sum\limits_{j=1}^k w_j \bm{A}_j\right)\right\Vert_{\ell}\geq\theta \right),
\end{eqnarray}
where $\left\Vert \cdot \right\Vert_{(\ell)}$ is Ky Fan $\ell$-norm. R.H.S. of Eq.~\eqref{eq2f h: UP sec: Function of Operator Bounds Algebra}, where the random tensors summation part is independent of functions $f$ and $h$ can be upper bounded by those theorems in Section IV in~\cite{chang2022generalMA}.

Consider the multiplication between $\sum\limits_{j=1}^k w_j \Phi(f(\bm{A}_j))$ and $\sum\limits_{j=1}^k w_j \Phi(h(\bm{A}_j))$, from Eq.~\eqref{eq2f: UP sec: Random Tensor Tail Bounds Algebra} and Eq.~\eqref{eq2h: UP sec: Random Tensor Tail Bounds Algebra} with assumptions of positive $\alpha_{f,\mathscr{U}}$ and $\alpha_{h,\mathscr{U}}$ and $g\left(\sum\limits_{j=1}^k w_j \bm{A}_j\right) > \bm{O}$ , then, we have
\begin{eqnarray}\label{eq3f h: UP sec: Random Tensor Tail Bounds Algebra} 
\mathrm{Pr}\left(\left\Vert\left(\sum\limits_{j=1}^k w_j \Phi(f(\bm{A}_j))\right)\times\left(\sum\limits_{j=1}^k w_j \Phi(h(\bm{A}_j))\right)\right\Vert_{\ell}\geq \theta\right)\nonumber \\
\leq \mathrm{Pr}\left(\left\Vert (\alpha_{f,\mathscr{U}} \times \alpha_{h,\mathscr{U}})\left(g\left(\sum\limits_{j=1}^k w_j \bm{A}_j\right)\right)^2\right\Vert_{\ell}\geq \theta\right)
\nonumber \\
=\mathrm{Pr}\left(\left\Vert g\left(\sum\limits_{j=1}^k w_j \bm{A}_j\right)\right\Vert_{\ell}\geq \sqrt{\frac{\theta}{(\alpha_{f,\mathscr{U}} \times \alpha_{h,\mathscr{U}})}}\right),
\end{eqnarray}
where the last term of Eq.~\eqref{eq3f h: UP sec: Random Tensor Tail Bounds Algebra} can be upper bounded by those theorems in Section IV in~\cite{chang2022generalMA}.

\bibliographystyle{IEEETran}
\bibliography{GenConv_Jensen_Approx_TailBoundsAlgebra_Bib}

\end{document}